\documentstyle[psfig]{article} %[10pt]
 \newcounter{remark}

\begin{document}
\newcommand{\unl}{\underline}
\newcommand{\8}{\partial}
\newcommand{\dis}{\displaystyle}
\newcommand{\beq}{\begin{equation}}
\newcommand{\eeq}{\end{equation}}
\newcommand{\Cset}{{ C}\!\!\!\!{I}\;}
\newcommand{\Rset}{{\rm I}\!{\rm R}}
\newcommand{\barr}{\begin{array}}
\newcommand{\earr}{\end{array}}
\newcommand{\remark}{\vspace{1mm} \noindent \stepcounter{remark}
{\bf Remark \arabic{remark}: }}

\title{Equation-free optimal switching policies for bistable reacting systems
using coarse time-steppers.}
\author{
 Antonios Armaou\thanks{Department of Chemical Engineering,
The Pennsylvania State University, State College, PA 16802,
E-mail:armaou@engr.psu.edu},
 Ioannis G. Kevrekidis\thanks{Department of Chemical Engineering,
PACM and Mathematics, Princeton University, Princeton, NJ 08544,
E-mail: yannis@princeton.edu}} \maketitle

\begin{abstract}
We present a computer-assisted approach to locating approximate
coarse optimal switching policies between stationary states of
chemically reacting systems described by microscopic/stochastic
evolution rules.
The ``coarse time-stepper" constitutes a bridge between the
underlying kinetic Monte Carlo simulation and traditional,
continuum numerical optimization techniques formulated in discrete
time.
The approach is illustrated through two simple catalytic surface
reaction models, implemented through kinetic Monte Carlo: $NO$
reduction on $Pt$, and  $CO$ oxidation on $Pt$.
The objective sought in both cases is to switch between two
coexisting stable stationary states by minimal manipulation of a
macroscopic system parameter.
\end{abstract}

\begin{center}
keywords: equation-free, dynamic optimization, coarse
timesteppers, kinetic Monte-Carlo, Hooke-Jeeves\\
\end{center}

\section{Introduction}

The computation of operating conditions that are
optimal with respect to a predetermined objective (derived from
economic considerations, safety requirements and product quality
constraints) has been, for many decades, an exciting research area
for chemically reacting systems.
Increase in computational power has led to increased attention to
dynamically evolving processes and the search for optimal
time-varying operation protocols, both for lumped and spatially
distributed systems (e.g.,
\cite{VVB90,vsp94,ac02,va_aiche04,cbav04}).
Interestingly, as sensing and actuation
become increasingly more resolved in space and time for emerging
chemical processes, spatiotemporally complicated operating
policies can be considered (e.g., \cite{WPKRE01,PWKRE02,Grier03}).
The problem of computing optimal transition paths over complicated
free energy surfaces (see e.g., \cite{bcdg02}) is also
computationally related, especially considering Jarzynski's
relation between equilibrium free energy and nonequilibrium work
(see e.g., \cite{sun03}).

Existing computational approaches for solving dynamic optimization
problems at the deterministic, continuum level may involve (a)
formulation as a temporally discretized problem (both for the
process and for the operating variable(s)) simultaneously, and
solution using large sparse linear algebra techniques (e.g
\cite{vasil93,bcw02,VVB90}); (b) formulations involving direct
integration of the model equations in time, keeping track of
possible constraint violations and temporal discretization of the
operating variables \cite{bcvm00,fb99,vsp94,cbav00}; or possibly
(c) using efficient solution algorithms within a dynamic
programming formulation \cite{ks00,Luus00,bertsekas87}.
Knowledge of a macroscopic process model, in the form
of macroscopic mass balances closed through appropriate constitutive
expressions -such as chemical kinetic rate formulas-,
is a fundamental prerequisite for these computational solution strategies.

In contemporary engineering applications, however, we are often
faced with problems for which the available physical description
is in the form of atomistic / stochastic evolution rules (e.g.,
kinetic Monte Carlo, Lattice Boltzmann, molecular dynamics,
Brownian dynamics) while the design, optimization or control is
required at a coarse-grained, macroscopic level.
The lack of a closed-form process model can in principle be
circumvented during the formulation of steady state optimization
problems through the construction of identification algorithms
``wrapped" around black box simulators (e.g., see \cite{spk03}).
Through design of computational experiments, sampling and
estimation, one can determine the form of ``surrogate'' models
that are subsequently used to identify local or global (e.g., see
\cite{jsw98,mf02}) optima.
In another approach, near optimal solutions are found for
unconstrained optimization problems, under the assumption that
absolute and relative error bounds are known for the computed
objective function values \cite{ks03}.

Over the last few years we have been developing an {\it equation-free} computational
approach enabling microscopic / stochastic simulators to {\it
directly} perform system-level tasks, such as coarse integration,
stability analysis, bifurcation/continuation and feedback control,
thus circumventing the derivation of explicit macroscopic
evolution equations \cite{TQK00, GKT02, KGHKRT02, SAMK02}.
In this work we demonstrate the applicability of this
computational enabling technology to coarse-grained optimization
tasks; in particular, we present a formulation methodology capable
of addressing dynamically evolving processes, to derive
optimization formulations that can be solved employing standard,
off-the-shelf direct search algorithms.
The approach is applied to computationally identify coarse optimal
operating parameter policies capable of switching the (expected)
behavior of kinetic Monte Carlo simulators from one stationary
state to another (alternatively, from the bottom of one potential
well to the bottom of another); this is illustrated  through
kinetic Monte Carlo (KMC) simulations of two simplified models of
heterogeneous catalytic chemical reactions.
Both these systems are  characterized by two ``coarse-grained" stable
stationary states.
We seek optimal (for a particular definition of
a cost function) parameter variation policies that will switch
the KMC simulation from one stationary state to the other within a
finite time interval.

The paper is organized as follows: In the following Section we present
the problem we want to solve and our two illustrative examples.
We then present our formulation for coarse-grained computational
optimization.
In the following two Sections we present detailed computational results
for each of our illustrative problems; we then conclude with a brief discussion.

\section{Process description}
We investigate microscopic/stochastic processes for which we
believe that the coarse-grained, expected dynamics can be
well approximated in closed form by an equation of the general type
 \beq\label{process}
\dot{x} = F(x,p)
 \eeq
but where the right-hand-side of the evolution law, $F$, is
not available in closed form.
Here, $x(t) \in\Rset^n$ is a state variable vector, $t$ is
the time, $\dot{x}$ is the time derivative of $x$,
$\dis{\frac{dx}{dt}}$, and $p \in {\cal P}^m$ is the vector of
process parameters (${\cal P}^m \subset \Rset^m$ is the subset of
accessible values of the process parameters).
The state variables $x(t)$ in the model of Eq.\ref{process} are
coarse-grained observables of the microscopic simulation -
typically a few lower order moments of an atomistically or
stochastically evolving distribution (e.g., a concentration, or,
in our heterogeneous catalytic reaction examples, a surface
coverage, the zeroth moment of the distribution of adsorbates on
the surface).
%
%{\bf the next sentence is a mess, and all of this is a mess -- a
%microscopic process is NOT at a steady state, it is "vibrating"
%around a stationary state. You have written the paper as if the
%expectation of the process and the process is the same thing. It
%steady state is for "the process" - let along that "the process"
%has NO steady states}.
% ANSWER: The classic problem when I've been hanging around models too long.
% Addition: Similarly, steady states of the process model
%           are representations of atomistic distributions around which
%           the process reaches a dynamic equilibrium
%
The (unavailable) equation for the expected behavior of the
process may possess, at fixed process parameter values, one or
more steady states $x_{ss,i}$.
We will assume that initially the microscopic process is in the
neighborhood of one of these stable, coarse-grained stationary
states, represented by $x_{ss,1}$
for the system of Eq.\ref{process}.

\subsection{The coarse time-stepper}

The basis of our approach is the computation of a deterministic
optimal policy -a time-dependent protocol for changing
one of the operating parameters of the process- so as to attain
a particular goal for the expected dynamics of the process (the
``coarse-grained dynamics"): switching from one stationary state
to another one.
We want to accomplish this by acting directly on the microscopic
simulator, thus circumventing the necessity of first deriving a
closed form evolution equation for these coarse-grained dynamics.
This deterministic policy for the coarse-grained behavior will
then be applied to individual realizations of the process.
As we will discuss below, it is convenient in our approach to reformulate
the coarse dynamics in discrete rather than continuous time.
The coarse-grained evolution law then takes the form:
 \beq\label{proc_disc}
 \barr{l}
x_{i+1} = G_T(x_{i},p_{i+1})  \vspace{2mm}\\
 \earr
 \eeq
where $t_i$ is the time at the beginning of $i$-th time interval,
$T$ is the time interval duration (reporting horizon), $G_T$
represents the evolution of Eq.\ref{process} for constant process
parameter value $p(t)=p_{i+1},\; \forall t\in(t_i,t_i+T]$,
initialized at $x_i$ and evolved for time $T$, arriving at state
$x_{i+1}$.

Conventional algorithms for the solution of optimization problems
involving $G_T$ incorporate frequent calls to a subroutine that {\it evaluates}
$G_T$ and/or the action of its derivatives in phase and parameter
space on initial conditions.
Equation-free based algorithms {\it estimate} the same quantities
by short bursts of (possibly ensembles of) microscopic simulations
{\it conditioned on} the same macroscopic initial conditions
(\cite {TQK00, GKT02, KGHKRT02, SAMK02, kgh04}).
The coarse time-stepper constitutes such an {\it estimate} of the
discrete-time, macroscopic input-output map $G_T$ obtained via the kinetic
Monte Carlo simulator.
Through a {\it lifting} operator the macroscopic initial condition
is translated into several consistent microscopic initial conditions
(distributions conditioned on a few of their lower moments).
In some cases constructing microscopic initial conditions
consistent with macroscopic observable quantities is easy (e.g.,
constructing distributions with prescribed means and variances);
in other cases (e.g., when we want to prescribe pair probabilities
on a lattice, or even a pair correlation function) an optimization
problem may need to be solved to successfully lift; short runs with
constrained dynamics algorithms
(\cite{SHAKE,KevGear2004a,KevGear2004b}) can also help in
preparing such consistent microscopic initial conditions.
%    The KevGear2004a  is     J. Sci. Comp.  in press

%
This ensemble of microscopic initial conditions is then evolved
microscopically, in an easily parallelizable fashion (one
consistent realization per CPU) for a relatively short time.
The results are averaged through a {\it restriction} operator back to
a macroscopic ``output"; it is precisely this output that traditional
algorithms simply compute through function evaluations when the
evolution equations are available in closed form.
As extensively discussed in \cite{mmk02,KGHKRT02} part of the
microscopic evolution is spent in a ``healing" process - the
higher moments, which have been initialized ``wrong" quickly relax
to functionals of the low order moments (our state variables).
A separation of time scales (fast relaxation of the high moments
to functionals of the low ones, and slow -deterministic- evolution
of the low ones) underpins the existence of deterministic
coarse grained evolution laws.
The dynamics of the evolving microscopic distribution
moments in the problems we study constitute thus a singularly perturbed system.
The requirement of finite time microscopic evolution (necessary for
the moment healing process) conforms with the discrete-time
formulation of the coarse optimization problem, which
is common in many optimization algorithms (see section \ref{s:opt}
below).

\subsection{Numerical experiments}\label{s:proc}

We illustrate the proposed combination of traditional optimization
techniques with stochastic simulators through KMC realizations
(using the stochastic simulation algorithm, proposed by Gillespie
\cite{G76,G77,G92}) of kinetic models describing catalytic $NO$
reduction by $H_2$ and $CO$ oxidation by $O_2$, respectively.

We initially apply the proposed method to a drastically simplified
kinetic model of $NO$ reduction by $H_2$ on Pt surfaces; the model
involves Langmuir adsorption, first order desorption, and chemical
reaction requiring two neighboring vacant sites.
The mean field Langmuir-Hinshelwood approximation of the kinetic mechanism
we will model consists of a single Ordinary Differential Equation for the
surface coverage of $NO$:
 \beq\label{rxn_ct}
 \dis{\frac{d\theta}{dt}}=\alpha(1-\theta)-\gamma
 \theta-k(1-\theta)^2 \theta.
 \eeq
Here $\theta$ describes the surface coverage of adsorbed $NO$,
$\alpha$, $\gamma$ are the $NO$ adsorption rate (incorporating
the dependence of gas phase pressure of $NO$) and desorption rate
constants respectively, and $k$ is the reaction rate
constant.
The reaction term is third order due to
the need for two free adjacent sites for
the adsorption of $H_2$.
In Figure \ref{fig1} we present the deterministic bifurcation
diagram in the form of coverage $\theta_{ss}$ at steady-state as a
function of $k$ for $\alpha=1$ and $\gamma=0.01$.
We observe a range of $k$ values for which
the system exhibits multiple steady-states; the higher and
lower ones are locally stable, while the middle one is
unstable.

The proposed approach is also applied to a kinetic Monte Carlo
description of the $CO$ oxidation reaction mechanism, using a simplified
model for the reaction kinetics of the form $A+1/2B_2 \rightarrow
AB$.
The kinetics here involve Langmuir adsorption for $A$, dissociative
adsorption for $B_2$, and a second order surface reaction whose products
are immediately desorbed.
The mean field approximation equations for this mechanism in the
absence of adsorbate interactions would consist of a set of two
ODEs \cite{mmk02}:
 \beq\label{rxn_2ode}
 \barr{r@{}l}
 \dis{\frac{d\theta_A}{dt}}=\alpha(1-\theta_A-\theta_B)-\gamma\theta_a-
        4k_r\theta_A\theta_B \vspace{2mm}\\
\dis{\frac{d\theta_B}{dt}}=2\beta(1-\theta_A-\theta_B)^2-4k_r\theta_A\theta_B
 \earr
 \eeq
where $\theta_A$ and $\theta_B$ describe the surface coverage of
adsorbed $CO$ and $O_2$ respectively, $\alpha$ and $\beta$ denote
the adsorption rate constants of $CO$ and $O_2$ respectively,
$\gamma$ the desorption rate constant of $CO$, and with $k_r$ we
denote the oxidation rate constant.
In Figure \ref{fig2} we present the deterministic bifurcation
diagram representing the coverage of A, $\theta_A$, and B,
$\theta_B$, respectively, as a function of the adsorption rate
constant $\beta$ for values of $\alpha=1.6$, $\gamma=0.04$ and
$k_r=1$.
We observe that the system exhibits multiple steady
states, the higher and lower ones being locally stable and the
middle one unstable.
Note that when the steady-state value of A
coverage is low, then the coverage of B is high and vice-versa.
%
%For the specific processes, we use Kinetic Monte-Carlo simulators
%to describe the dynamic behavior of the process. Specifically, we
%assume that we know the mean field evolution equations for the
%coverage $\theta_i$ of each species $i$ on the catalytic surface
%\cite{mmk02}.
%
A description of the use of the coarse KMC time-stepper in
obtaining ``coarse" versions of this bifurcation diagram may be
found in \cite{mmk02}.

%-------
\section{Coarse computational optimization}\label{s:opt}

Computing the temporal profiles of the process parameters that
cause the transition of the coarse-grained system from one initial
stationary state to a different final one can be formulated as an
optimization problem; the objective is to minimize an integral
cost function over time:
 \beq\label{opt_obj_inf}
 \dis{\min_{p(t)\in {\cal P}^m}}\dis{\int_0^{\infty}}
{\cal Q}(t,x,p)
 dt
 \eeq
where ${\cal Q}$ is a continuous scalar cost function. The
constraints for this coarse optimization problem are the
(unavailable) coarse process evolution equations Eq.\ref{process},
the initialization at one of the coarse stationary states for $p =
p_{ss}$, the requirement of termination at a different stationary
state for the same process parameter value $p_{ss}$, as well as,
possibly, other inequality constraints $g(x,p)$ (e.g., that
surface coverages should remain positive and sum up to less than 1
at all times):
 \beq\label{opt_con_inf}
 \barr{c}
 \dot{x}-F(x(t),p(t))=0, \quad g(x,p)\le 0 \vspace{2mm}\\
 \barr{ll}
 p(0)=p_{ss},& \dis{\lim_{t\rightarrow \infty}}p(t)=p_{ss} \vspace{2mm}\\
 x(0)=x_{ss,1},& \dis{\lim_{t\rightarrow \infty}}x(t)=x_{ss,2}.
 \earr
 \earr
 \eeq
This constitutes an infinite dimensional problem in continuous time.
Direct
solution methods are based on the calculus of variations.
Semi-infinite programming approaches provide us with the necessary
mathematical tools to solve such problems with finite time horizon
\cite{hk93}, through discretization of the temporal domain.

Another approach consists of approximating this problem through
a finite time horizon problem with a final state penalty,
which is (in our case) solved in discrete time.
This results in a finite dimensional,
generally nonlinear, optimization problem which
-if the coarse equations are available-
could be solved using available optimization techniques.
For example, discretizing the process time $[0,t_f]$ in $N$ time
intervals of length $T$ (not necessarily constant) and assuming
that the process parameters remain constant within each interval,
results in the following optimization program with
$(N+1)\times(n+m)$ variables and $n\times (N+1)+m+n$ equality
constraints:
 \beq\label{opt_fin-msm}
 \barr{c}
\barr{r@{}l} \dis{\min_{p_i\in {\cal P}^m}}& \dis{\sum_{i=1}^N}
\dis{\int_{(i-1)T}^{iT}}{\cal Q}_d(t,x_i,x_{i-1},p_i)dt \vspace{2mm}\\
&+{\cal W}(R(|x(t_f)-x_{ss,2}|-\epsilon))\earr\vspace{2mm}\\
 {\rm s.t.} \vspace{2mm}\\
 x_{i} = G_T(x_{i-1},p_i)=0,\quad i=1,2,\dots,N\vspace{2mm}\\
 g_d(x_{0},\dots,x_{N},p_0,\dots,p_N)\le 0,\vspace{2mm}\\
p_0=p_{ss},\quad  x_0=x_{ss,1}.
 \earr
 \eeq
Here ${\cal Q}_d$ is analogous to the ${\cal Q}$ function for the
discrete values of the state, $g_d$ is similarly analogous to $g$,
${\cal W}(\cdot)$ is a class $K$ scalar function, $R(\cdot)$ is
the ramp function and $\epsilon$ is a value for which
$|x_N-x_{ss,2}| \le \epsilon <=> {\cal W}=0$.
Appropriate final state penalty contributions to the
objective function take the place of final state constraints at
infinite time as stated in the original formulation of the
problem; the final state is restricted to be (in finite time)
within a neighborhood of the final stationary state.
This fully discrete-time formulation is ideally suited for linking
with a coarse time-stepper.
It is also in principle possible to discretize only the process
parameter temporal behavior, resulting in the following
formulation (coined {\it control vector parameterization}
\cite{vasil93}) with $N\times m+2n$ variables and $n$ equality
constraints:
 \beq\label{opt_fin-cvp}
 \barr{c}
\barr{r@{}l} \dis{\min_{p\in {\cal P}^m}}&
\dis{\int_{0}^{t_f}}{\cal Q}(t,x,p)dt \vspace{2mm}\\
&+{\cal W}(R(|x(t_f)-x_{ss,2}|-\epsilon))\earr\vspace{2mm}\\
 {\rm s.t.} \vspace{2mm}\\
\dot{x}-F(x,p)=0,\quad g(x,p)\le 0
 \vspace{2mm}\\
% p(t)=\dis{\sum_{i=1}^{N}p_iH(\frac{T}{2}-|t-(i-\frac{1}{2})T|)}+p_{ss}H(-t)
p(t)=\dis{\sum_{i=1}^{N}p_i\Pi(\frac{t}{T}-i+\frac{1}{2})}+p_{ss}H(-t)
 \earr
 \eeq
where $\Pi(\cdot)$ is the standard boxcar function and $H(\cdot)$
denotes the Heaviside function.
In cases where the explicit form of Eq.\ref{process} is
unavailable, the state evolution is provided through direct
simulation of the system.

\subsection{Solution methodology}\label{s:sol}

Traditional discrete time optimization schemes would repeatedly call, during
the solution process, a numerical integration subroutine for the
system equations (and possibly variational and sensitivity
integrations for the
estimation of derivatives with respect to state variables and
parameters).
This call is now substituted by the coarse time-stepper; the most
important numerical issue is that of noise, inherent in the
lifting process and the stochastic simulations, and the variance
reduction necessary to estimate the state or its various
derivatives.
Simulations of different physical size systems (different lattice
sizes in our simulation) are characterized by different values of
variance.
For the type of simulations in this paper, changing the physical
size of the simulated domain on the one hand, and changing the
number of copies of the simulation on the other, have comparable
effects in reducing the variance of the simulation output; this,
however, is not generally the case in KMC simulations.
When a switching policy {\it for a particular physical size
system} is required (e.g., for nanoscopic reacting systems, such
as the chemical oscillations on Field Emitter tips \cite{SBJEI99})
variance reduction can be affected through a larger ensemble of
consistent microscopic initializations (see also \cite{MO95,O96}
for variance reduction techniques).
Simulation noise affects both function evaluations and
(numerical) coarse derivative evaluations, and thus becomes
an important element of the approach.
It is interesting, however, to observe that massively parallel computation
can reduce the {\it wall clock time} required for the computation
distributing different microscopic initializations/realizations to
different CPUs.

Due to the presence of noise in the coarse time-stepper results,
the use of optimization algorithms that are specifically designed
to be insensitive to noise becomes desirable and even necessary
(the numerical estimation of derivatives is, of course, highly susceptible
to noise).
A class of search algorithms that fulfill this requirement are ones
that use only function evaluations to search for the optimum, such
as the Luus-Jaakola \cite{lj73}, Nelder-Mead and Hooke-Jeeves
algorithms.
Algorithms that use local and bounded approximations of the
Jacobian matrix of the objective function with respect to the
process parameters have also been developed, such as the implicit
filtering algorithm; the reader may refer to \cite{kelley99,klt03}
for a review of these methods.
In our approach, an optimization ``wrapper" is constructed around
an {\it on demand estimate} of the system coarse-grained model, constructed
by computational experimentation with the microscopic solver.
There is a close relation to the algorithms in
\cite{kelley99,klt03,ks03}, where comparable wrappers are
constructed around {\it continuum} models that are difficult or
expensive to evaluate (e.g., the ones that come from solutions of
large PDEs).
It is also remarkable that this ``wrapper" approach, which performs
the equation-free solution of an optimization problem through appropriate
design of {\it computational} experiments, can be also in principle be
used for the
equation-free optimization of {\it experimental} systems.
Indeed, if enough control authority exists for experiments to be
initialized in detail, these optimization ``wrappers" become
protocols for the design of laboratory experiments leading to an
optimum \cite{gbbb97,lmr01,rhr04}.

In all the above iterative algorithms \cite{kelley99}, an
appropriate initial guess of the process policy profile and a set
of search directions for the $m\times N$ coarse-grained
variables is required.
The usual set of search directions (also used in our numerical
experiment) is the unitary basis for $\Rset^{m\times N}$.
Specifically for stencil-based search algorithms, important
parameters also include a vector, the elements of which are the
maximum distances the algorithm should venture from the current
position during the new direction search step at each iteration.
These perturbation distances are called {\it scales}.

Selection of these scales in the search algorithm requires estimates
of noise bounds.
A scale that is ``too small" not only leads to
increased computation time with no apparent advantages, but,
especially in the case of algorithms that compute approximations
of the Jacobian, can lead to grossly erroneous results for the
search.
The following time-stepper protocol yields, in our case,
simulation results of bounded ensemble variance:
\begin{enumerate}
 \item Initialize time-stepper. Set:
 \begin{itemize}
  \item lattice size $N_{l}$,
  \item solution replica ensemble size $M_{r}$,
  \item minimum $M_{min}$ and maximum $M_{max}$. % $N_{max}$ and
  \item ensemble variance limit $d_{max}$ and

%%%YANNIS
% Antoni, noise amplitude is too loose a word, we need something
%  more careful--- something like variance of the ensemble results
% or something -- "noise amplitude" is funny
%%%%%%%%%%%%%%%%%%%
%- We could call it variance limit. The thing is that the way we have it
%- (which I am not certain it is the best) is variance of each single variate.
%- This of course is not the variance of the multiple variables we may have.
%- How ever, I would think it is the best to consider the variance of the multivariate
%- system. BY THE WAY, KELLEY AND SACHS PUBLISHED THE SAME THING IN 2003
%- of course with a theoretical analysis.

 \end{itemize}
 \item In {\it each} time step, the time-stepper:
 \begin{enumerate}
  \item simulates system for the desired parameter value $M_{r}$ times (this can
    be affected through different microscopic initializations and/or
    different random seeds in the Monte Carlo process).
  \item Computes the standard deviation divided by the average of the sample, $d$.  \item If $d>d_{max}$ and $M_{r} < M_{max}$, adjusts $M_{r} \ge M_{min}$, repeat step (a).
 \end{enumerate}
\end{enumerate}
%
%- By the way, what the algorithm does is to remember how many copies it used
%- in the previous search for each timestep and use the relation between expected
%- variance and number of copies to adjust the initial guess of how many copies
%- are needed. If more are needed, it will recompute and increase accordingly.
%- We rewrote that part an fudged it because you were not convinced that
%- a relationship would always be available.
%- It was not blindly increasing the number of copies when it was needed.
%
Adaptively adjusting the number of realization in the ensemble
can thus help in the choice of appropriate scales.
The lattice size $N_l$ of the time-stepper routine affects the
expected evolution of the process observables; for the
SSA simulations in this paper, as the size of the simulation increases
the expected behavior of the stochastic microscopic
process becomes closer and closer approximated
by the mean-field description represented by the ODE system \cite{G77}.
In our computations in this paper, the lattice size of our
time-stepper is chosen so that the discrepancy between the
expected behavior obtained from the KMC simulations and the ODE
solutions is negligible \cite{fw91}; this way we can validate the
directly computed optimal switching policies by comparing them to
the mean-field, deterministic ones.
%
%DONE
%%%YANNIS -- this should be also better said, Antoni
%
%   if this size is large enough, the behavior of the
%   stochastic process approximates that of the mean-field ODEs
%  and the reference should be Gillespie, some of the papers where
%    he talks about the "Chemical Langevin" equation.
%%%%%%%%%%%%%%%%%%%
%
%Increase in lattice size has the main effect
%of altering the expected value of the process parameters, which
%converges to a certain value at the limit.
%
Changing the lattice size will, in general, affect the expected behavior
of a stochastic process; the variance of the time-stepper results will
also be affected.
We will use multiple replica simulation runs $M_r$  (with the
same process settings and macroscopic initial conditions) to control
fluctuations in the estimation of the expected results.
%attempt to control the variance in the obtained from the
%coarse timestepper results by manipulating the number of replica
%simulation runs $M_r$ of the process with the same process
%conditions and macroscopic initial conditions.

\subsection{Direct search method: The Hooke-Jeeves algorithm}\label{s:HJ}
The processes that we investigate often lead to optimization problems,
where the objective functions exhibit a large number of small
local optima that in effect ``hide" the directions along which the
objective function decreases.
This is due to the stochastic nature
of the microscale simulations that are used to infer the dynamic
behavior of the coarse variables.

A class of search algorithms that avails itself to the solution of
such problems is direct search algorithms.
These methods evaluate
the objective function at sample points and use the provided
information to continue sampling along promising directions.
A prerequisite is an initial point inside the feasible region as
well as a set of search directions.
A number of such methods exist
including Nelder-Mead, Hooke-Jeeves, Multi-Coordinate Search
(MCS); in the current work we use the Hooke-Jeeves method.
For
completeness we briefly describe the method (the reader is
referred to \cite{kelley99} for a detailed presentation).

Hooke-Jeeves is a stencil based method, similar to coordinate
descent; however it is a more aggressive search.
As in all stencil-based methods, a set of search directions
$\textbf{v}$ is provided, as well as an array of $M$ scales,
$\textbf{s}$, that determine the step length size which the search
algorithm is allowed to venture away from the current point when
investigating for promising search directions.
Let us define the design variable vector as $x\in\Rset^N$, the
objective function as $f(\cdot)$, and (without loss of generality)
the optimization problem as one of minimizing the objective
function.
In the following pseudocode we present the Hooke-Jeeves method,
comprising of the following steps \cite{kelley99}:

\begin{itemize}
 \item Define: Initial position, $x_0$, search directions, $\textbf{v}$, and scales, $\bf{s}$.
 \item search:
\begin{enumerate}
 \item Compute $f(x_0)$.
 \item choose scale $s_i$ from $\textbf{s}$.
 \item Exploratory step:
 \begin{enumerate}
 \item[3.1] Define $x_s=x_0$.
 \item[3.2] Investigate search direction $v_j\in \textbf{v}$:
 \begin{itemize}
 \item Compute $f(x_s+s_iv_j)$. If $f(x_s+s_iv_j)<f(x_0)$ define
 $x_s=x_s+s_iv_j$, move to step 3.3
 \item If $f(x_s+s_iv_j)\ge f(x_0)$ compute $f(x_0-s_iv_j)$.
 \item If $f(x_s-s_iv_j)<f(x_0)$ define
 $x_s=x_s-s_iv_j$, continue to step 3.3.
 \end{itemize}
 \item[3.3] Repeat step 3.2 for $j=j+1$. If $j=N$ (all $N$ search directions in $\textbf{v}$ investigated) continue to step 3.4.
 \item[3.4] Obtain promising direction $d_i=x_s-x_0$.
 \end{enumerate}
 \item If $d_i=0$, then $i=i+1$ (reduce scale size to $s_{i+1}$). If $i=M$
 (all $M$ scales in $\textbf{s}$ investigated) terminate search iterations.
 \item Pattern move step:
 \begin{itemize}
 \item Define $x_c=x_0+2d_i$. Compute $f(x_c)$.
 \item If $f(x_c)< f(x_s)$, set $x_s=x_c$.
 \end{itemize}
 \item Set $x_0=x_s$, and repeat step 1.
\end{enumerate}
 \item Optimal point, $x_0$, obtained.
 \end{itemize}

\begin{figure}[htbp]
\centerline{\psfig{file=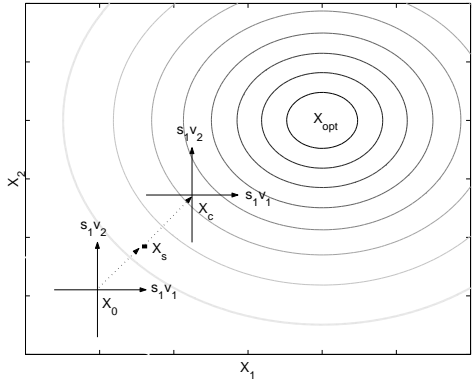,height=2.0in}}%,width=1.2in
\caption{Two-dimensional optimization problem, presenting the
steps of one iteration during an optimal search using Hooke-Jeeves
algorithm. Darker contour lines denote lower values of objective
function.} \label{figHJ}
\end{figure}

\begin{figure}[htbp]
\centerline{\psfig{file=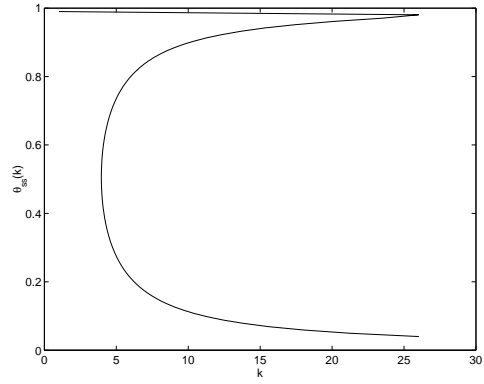,height=2.0in}}%,width=1.2in
\caption{Bifurcation diagram of $\theta$ at steady state with
respect to $k$.} \label{fig1}
\end{figure}

\begin{figure}[htbp]
\centerline{\psfig{file=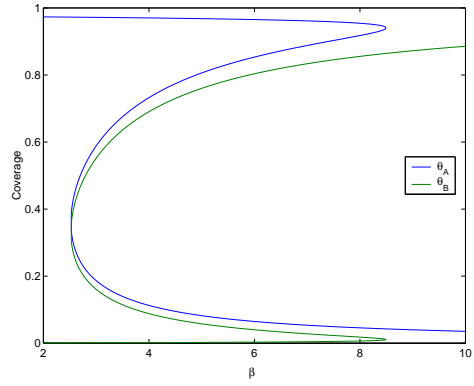,height=2.0in}}%,width=1.2in
\caption{Bifurcation diagram of $\theta_A$ and $\theta_B$ at
steady state with respect to $\beta$.} \label{fig2}
\end{figure}

In Figure \ref{figHJ} we sketch the application of the above
pseudocode in a representative optimization problem with $N=2$;
only one iteration is presented.
The optimal location in Figure \ref{figHJ} is represented by
$x_{opt}$, and darker contour lines present decreasing values of
the objective function.
Position $x_0$ is chosen to initialize the search (step 1), and
during the exploratory step (step 3), the objective function is
evaluated in the search directions $v_1$ and $v_2$ and at a
distance $s_1$ from $x_0$.
Two directions are computed to lead to lower objective function
values (shown with arrows, step 3.2), thus defining $x_s$ at the
end of step 3.3 to be at the upper right, and proposing the
move direction shown with the dotted arrow in Figure
\ref{figHJ}.
During the pattern move step (step 5), an aggressive move to $x_c$
is proposed (shown in Figure \ref{figHJ}), and $f(x_c)$ is
computed to be less than $f(x_s)$. As a result, the next iteration
takes place at position $x_c$, shown in Figure \ref{figHJ}

The method can be easily modified to incorporate inequality
constraints for the variables $x$, by modifying the objective
function to increase when $x$ moves outside a feasible set.
The search algorithm behavior at the boundaries is reminiscent of
interior point methods.
Optimization algorithms that are based only on function evaluations
are not guaranteed to converge to a global minimum, or even to a
minimum.
This is due to the fact that the necessary optimality conditions
in the neighborhood of the result are {\it not computed}, and also
because a search direction may have been neglected.

Once the optimization algorithm has converged and produced an
``apparent best" operating policy profile, it may be prudent to
restart using the result of the previous search as an initial
guess.
%
%Restarting the search at ``large" scales causes a large
%perturbation, which may assist in finding a better local optimum.
{Restarting the search, which invokes again the "large" scales,
will cause large perturbations in the operating policy profile,
which may assist in finding a better local optimum.}
Once two consecutive runs have produced comparable results, the
free variable profile is declared ``optimal over all scales"
\cite{kelley99}.
We note that due to the use of KMC simulations, two simulation
runs will not produce the same value of the objective function for
the {\it same} policy.
This results in small variations of the identified policy profiles
in two consecutive searches, if small scales are used (scales that
lead to perturbations in the system state evolution that are
comparable to the noise of the KMC simulations).

%------------
\section{Numerical Results: NO reduction on Pt surface}\label{s:NOresults}

The approach outlined in section \ref{s:sol} was initially applied
towards the computation of an optimal switching policy (between
different stationary states) in our simplified $NO$ reduction
model.
We define a (slightly arbitrary, but useful for illustration
purposes) objective function in Eq.\ref{opt_fin-cvp} as:
 \[
\barr{l}
 {\cal Q}=(k(t)-k_{ss})^2(1-0.3e^{-t})T(\dis{\sum_{i=0}^{N}}\delta(t-iT))\vspace{2mm}\\
 {\cal W}=50[1-exp(-R(|\theta(t_N)-\theta_{ss,f}|-\epsilon))]

\earr
 \]
where $\epsilon=0.05$, $R$ denotes the ramp function  and $\delta$
denotes the standard Dirac function. We assume that the single
``manipulated variable" for this problem (the process parameter
$p$ of Eq.\ref{process}) is the (macroscopic) reaction rate constant $k$;
this affects the reaction probability in the microscopic simulation.
In continuum models the reaction rate constant could be manipulated, for
example, by changes in temperature, through its Arrhenius temperature
dependence; our adsorption rate constants could be varied by
changes in the gas phase pressure of a species.
The options that were used
for the specific optimization problem are presented in Table
\ref{t:proc}.
The form of ${\cal Q}$ reflects the formulation of the time-dependent
policy through a finite number of decisions, placing a penalty
at the time of decision.
${\cal W}$ is a class K function
penalizing the deviation of the final state from the desired
steady-state heavily initially, reaching a plateau at the maximum
deviation of $\theta$.

KMC simulations based on the Gillespie SSA algorithm formed the
basis of the coarse time-stepper that was used to estimate the
coarse-grained system response.
A variety of lattice and sample sizes
were used in estimating the dynamic behavior of the system.
Hooke-Jeeves was the algorithm of choice in our search for the optimal profile.
In Table \ref{t:HJ.MC-obj} we present the value of the objective
function for the computed optimal parameter profiles, through
which the effect of the error in the computed optimal profile is
implicitly quantified.
The lattice size $N_{l}$ was increased with no appreciable change
in the expected system response.

\begin{figure}[htb]
\centerline{a)\psfig{file=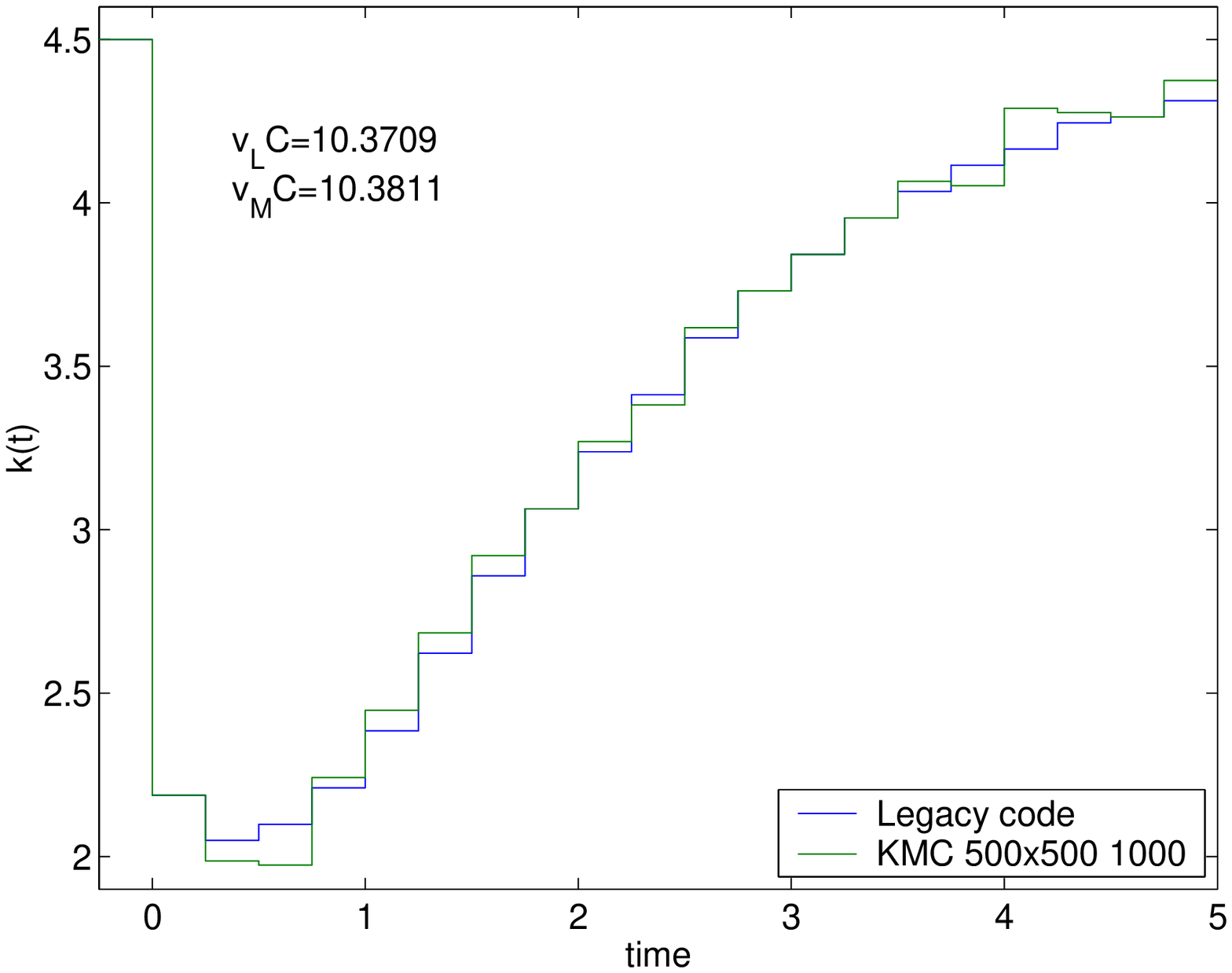,height=2.0in}\hspace{3mm}
b)\psfig{file=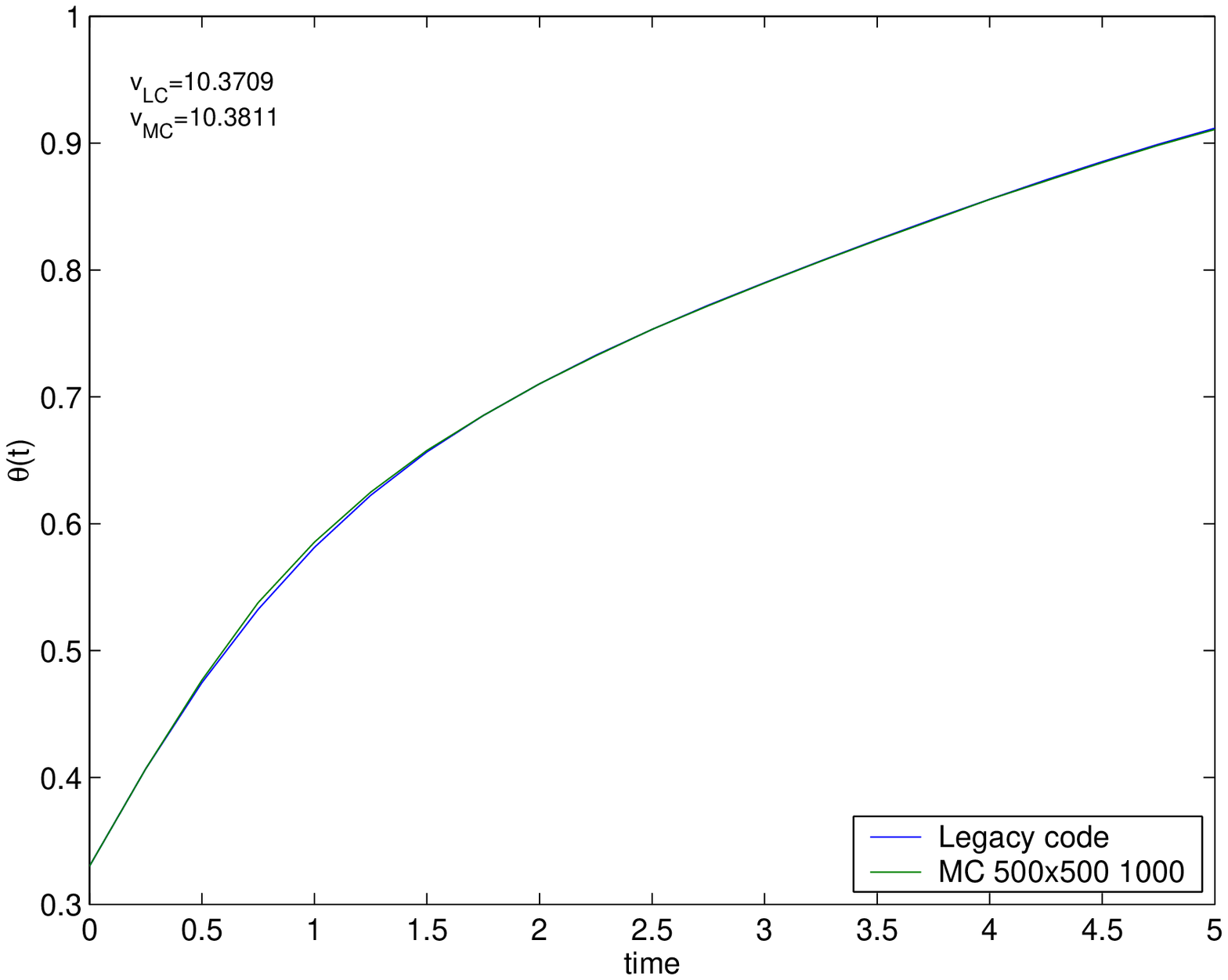,height=2.0in}}%,width=1.2in
\caption{Results using Hooke-Jeeves algorithm through numerical
integration of Eq.\ref{process} (blue line) and using KMC
simulation (green line), a) Optimal temporal profile of process
reaction rate $k$, b) Evolution of $NO$ coverage $\theta$.}
\label{figHJ2}
\end{figure}

A secondary gain from the lattice size increase was that the
ensemble variance $d$ decreased, since $ d \propto
(N_{l})^{-1/2}$.
Moreover, as the sample size $M_{r}$ used to compute the coarse
response increased in the KMC simulations, the ensemble variance
$d$ also decreased, as expected, since $d \propto (M_{r})^{-1/2}$.
This in turn lead to results from the search for the optimal
profile of $k$ that are closer to ones obtained when we use the
time-stepper of the actual deterministic problem (for comparison
purposes).
The Gillespie stochastic simulation algorithm \cite{G76,G92}
was chosen so that the coarse behavior at large sample sizes can
be well approximated through deterministic ODEs,
and so that the noisy time-stepper optimization results can be
compared to the deterministic ones at the appropriate limit.
The objective value convergence to the computed optimal value from
the ODE ``direct simulator" comes, of course, at the cost of increased
computational work.
This is only an illustrative example, to validate the approach; it
would not make sense to use coarse-grained computations when good
deterministic equations for the macroscopic behavior are known.
The equation-free approach is intended for cases where the macroscopic
equations are {\it not} available in closed form.
The use of parallel computing can, as we
discussed, drastically decrease the necessary wall-clock simulation time.
In Figure \ref{figHJ2} we present the
results for $N_{l}=500 \times 500$ and $M_{r}=1000$
and compare them to the direct macroscopic ODE simulation results.
A {\it near}-optimal parameter profile is arrived at, due to the
combination of KMC simulations' noise and the Hooke-Jeeves search
algorithm (a search direction that has been investigated and
characterized as unfavorable is not reinvestigated to conserve CPU
time).

We also used an alternative algorithm (Implicit Filtering,
\cite{kelley99}) to compute a near optimal steady state switching
profile of $k(t)$.
The method of implicit filtering uses consecutive bounded
approximations of the Jacobian employing a set (varying at each
iteration) perturbation of the design variables and poses a
prescribed limit to the maximum change of the design variable at
each iteration.
The time-stepper protocol used had $d_{max}=0.005$, with variable
lattice size set at $N_{l}=126\times 126$ and $N_{max}=4N_{l}$.
Adaptively adjusting $M_r$, originally set at $M_{r}=252$, with
$M_{max}=4M_{r}$ and $M_{min}=0.5M_r$, the bounds on the ensemble
variance were satisfied.
We also note that a minimum lattice size bound was always enforced
in order to ascertain that the discrepancy between the expectation
of the finite-size KMC simulations and the mean-field ODEs for the
same mechanism is negligible over the time scales we are working.

\begin{figure}[htb]
\centerline{a)\psfig{file=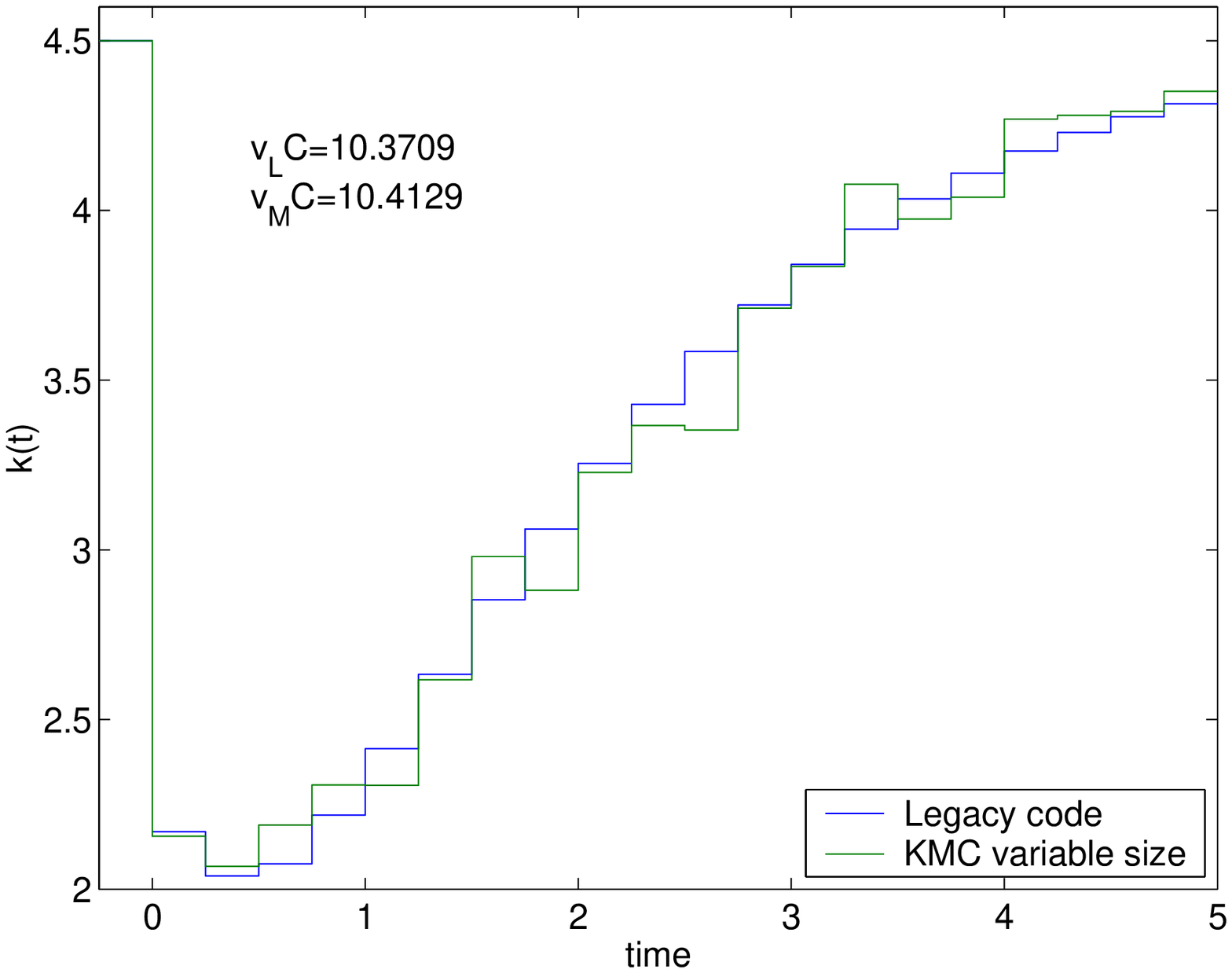,height=2.0in}\hspace{3mm}
b)\psfig{file=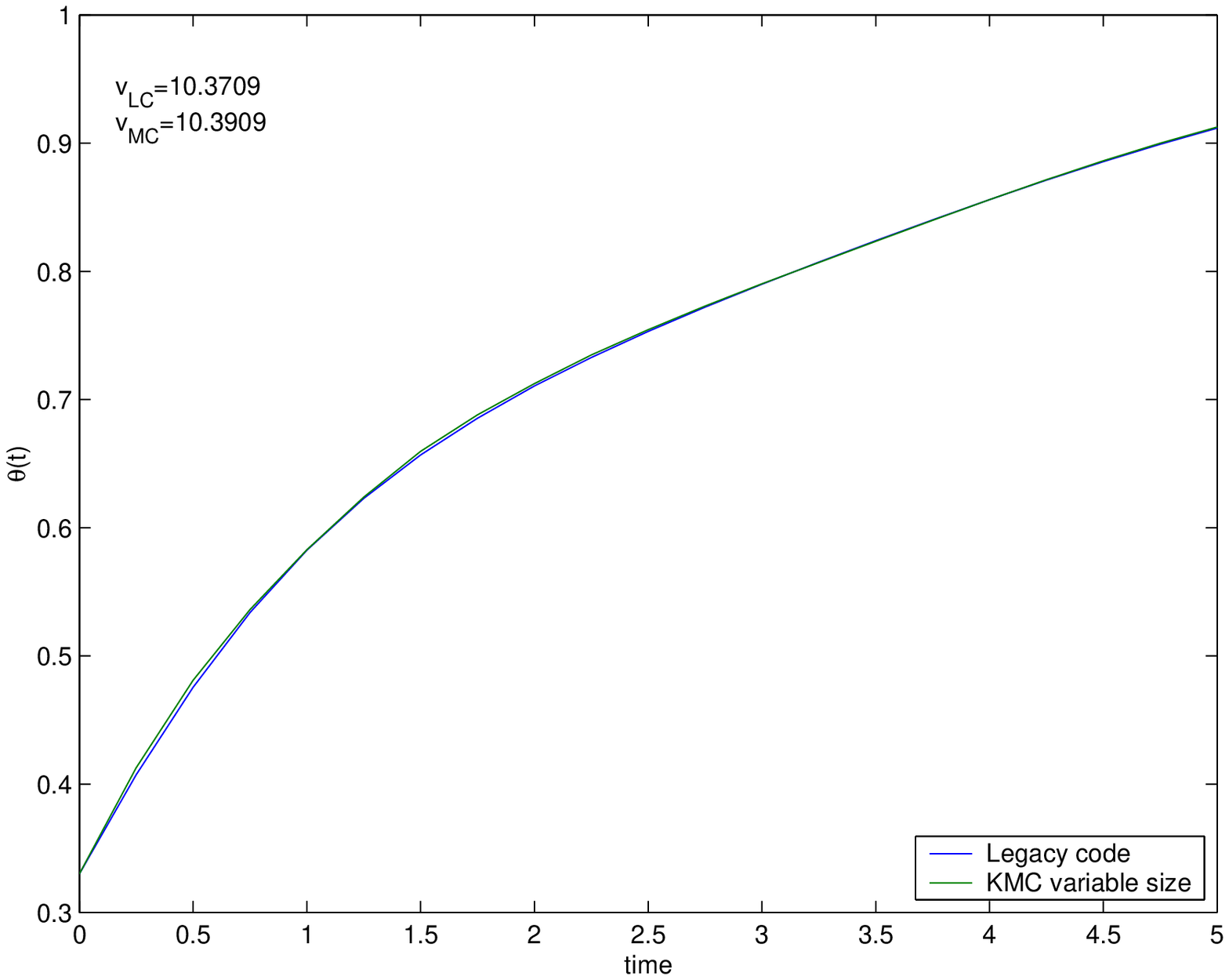,height=2.0in}}%,width=1.2in
\centerline{c)\psfig{file=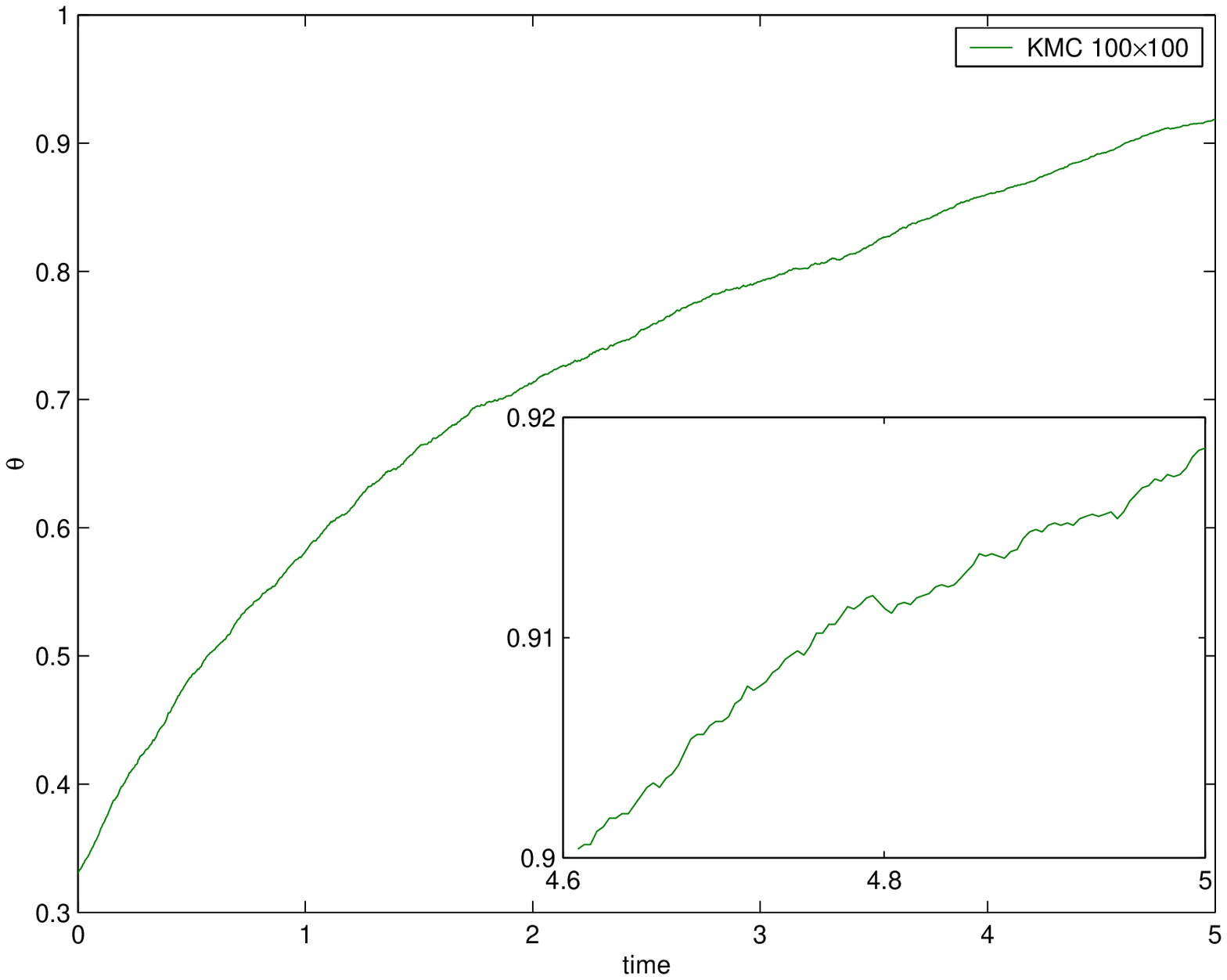,height=2.0in}}
\caption{Results using Implicit Filtering algorithm through
numerical integration of Eq.\ref{process} (blue line) and using
KMC simulation (green line), a) Optimal temporal profile of
process reaction rate $k$, b) Evolution of $NO$ coverage $\theta$,
c) Evolution of $NO$ coverage $\theta$ for a single KMC
realization, $N_l=100\times 100$, $\delta t=0.0039$. }
\label{figIF1}
\end{figure}

The resulting near optimal time profile of $k(t)$ is presented in
Figure \ref{figIF1}a, while the near optimal path of $NO$ coverage
evolution $\theta(t)$ is shown in Figure \ref{figIF1}b for an
averaged realization, and in Figure \ref{figIF1}c for a single KMC
realization with lattice size $N_l=100\times 100$ and reporting
horizon $\delta t=0.0039$.

In Table \ref{t:IF.MC-obj} we present computational results
obtained through incorporating the coarse time-stepper in an
implicit filtering algorithm.
As one might expect, as the variance limit is decreased the
algorithm -for a large enough lattice- approaches the optimal
profile of $k(t)$ computed directly through the coarse ODE model.
We also observe that a successful scales selection depends heavily
on $d_{max}$, as scales below a certain limit have an adverse
effect on the computed near optimal path result (compare the
results of the search when the lowest scale is $2^{-4}$ to the one
when it is $2^{-3}$).

The time-stepper protocol can be combined with the search
algorithm to solve the optimization program using successively
more refined scales.
A first search, with large values for the scales, can take place
at higher $d_{max}$, followed by searches with gradually lower
$d_{max}$ and smaller scales to refine the search for the optimal
path; this approach may lead to computational savings.
When using the KMC with $d_{max}=0.005$ and a lower scale of
$2^{-3}$ for an initial search, followed by a KMC simulation of
the system with $d_{max}=0.0005$ and lower scale of $2^{-5}$, we
find a near optimal path with $v=10.4110$ in a total time of
$3981\;s$.
A search using KMC with $d_{max}=0.001$ and lower scale of
$2^{-4}$ lead to computing a near optimal path of $v=10.4129$ in
$6858\; s$ (results are shown in Table \ref{t:IF.MC-obj}).

The ``optimal" switching path, shown in figure \ref{figIF1}b, takes
the (expected) phase point from $\theta_{ss,s}$ through the
unstable steady state $\theta_{ss,i}$ as shown in Figure
\ref{figIF1}b.
After this is accomplished, we observe that the
optimal $k(t)$ trajectory rapidly converges back to $k_{ss}$
(see Figure \ref{figIF1}a).
The coarse phase point is now within the region of attraction
of the steady state  $\theta_{ss,f}$, and no particular
switching action is needed to get it there.

The objective function values presented in this section were
computed, for comparison purposes, by integrating the coarse
system Eq.\ref{rxn_ct} for the optimal switching profile found.
During the optimal search using KMC simulations, the coarse
process model was {\it never} used.
Several optimization algorithms were explored in conjunction with
the coarse time-stepper, including Hooke-Jeeves, Nelder-Mead,
Implicit-Filtering as well as Multilevel Coordinate Search (MCS).
Hooke-Jeeves was primarily chosen due to the
simplicity of the method and its relative convergence speed.
Implicit Filtering algorithm was mainly used with the coarse
time-stepper to illustrate the variance-reduction protocol.
Enforcing an upper bound on the variance of the ensemble, along
with the selection of the value of the scales, provided an order of
magnitude estimate of the error in the estimated derivatives
during the Jacobian approximation; variance reduction techniques
may also be useful for this purpose.

%------------
\section{Numerical Results: catalytic CO oxidation}\label{s:COresults}

\begin{figure}[htbp]
\centerline{\psfig{file=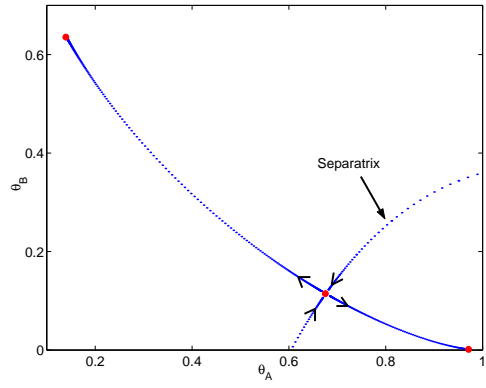,height=2.0in}}%,width=1.2in
\caption{phase portrait of system of Eq.\ref{rxn_2ode} presenting
the three steady states and the separatrix for $\beta=3.5$
(process parameter values shown in Table \ref{t:COproc}).}
\label{fig.2D.Mnfld}
\end{figure}

After solving a one-dimensional (coarse grained) example,  we
used the same computational approach
to find optimal switching policies between different
stationary states of a simplified $CO$ oxidation model presented
in section \ref{s:proc}.
The (again slightly arbitrary, illustrative)
objective function in Eq.\ref{opt_fin-cvp} was defined as:
 \[
\barr{l}
 {\cal Q}=(\beta(t)-\beta_{ss})^2(1-0.3e^{-t})T(\dis{\sum_{i=0}^{N}}\delta(t-iT))\vspace{2mm}\\
 {\cal
 W}=50[1-exp(-R(|\theta_{CO}(t_N)-\theta_{CO,ss,f}|-\epsilon)-R(|\theta_{O_2}(t_N)-\theta_{O_2,ss,f}|-\epsilon))]
\earr
 \]
where $\epsilon=0.05$, $R$ denotes the ramp function and $\delta$
denotes the Dirac function.
The single ``manipulated variable" for this problem
(the process parameter $p$ of Eq.\ref{process}) was chosen to
be the oxygen
$O_2$ adsorption rate constant $\beta$, which can in principle be
manipulated experimentally by varying the oxygen gas phase pressure.
The options that were used for the specific optimization problem
are presented in Table \ref{t:COproc}.
In figure \ref{fig.2D.Mnfld} we present the phase portrait of the
system under the specific nominal process parameter settings.
The evolving coarse-grained system
has to traverse a one-dimensional separatrix in two phase-space
dimensions, in order to enter the basin of attraction of the
desired $CO$ rich steady-state.

\begin{figure}[htbp]
\centerline{a)\psfig{file=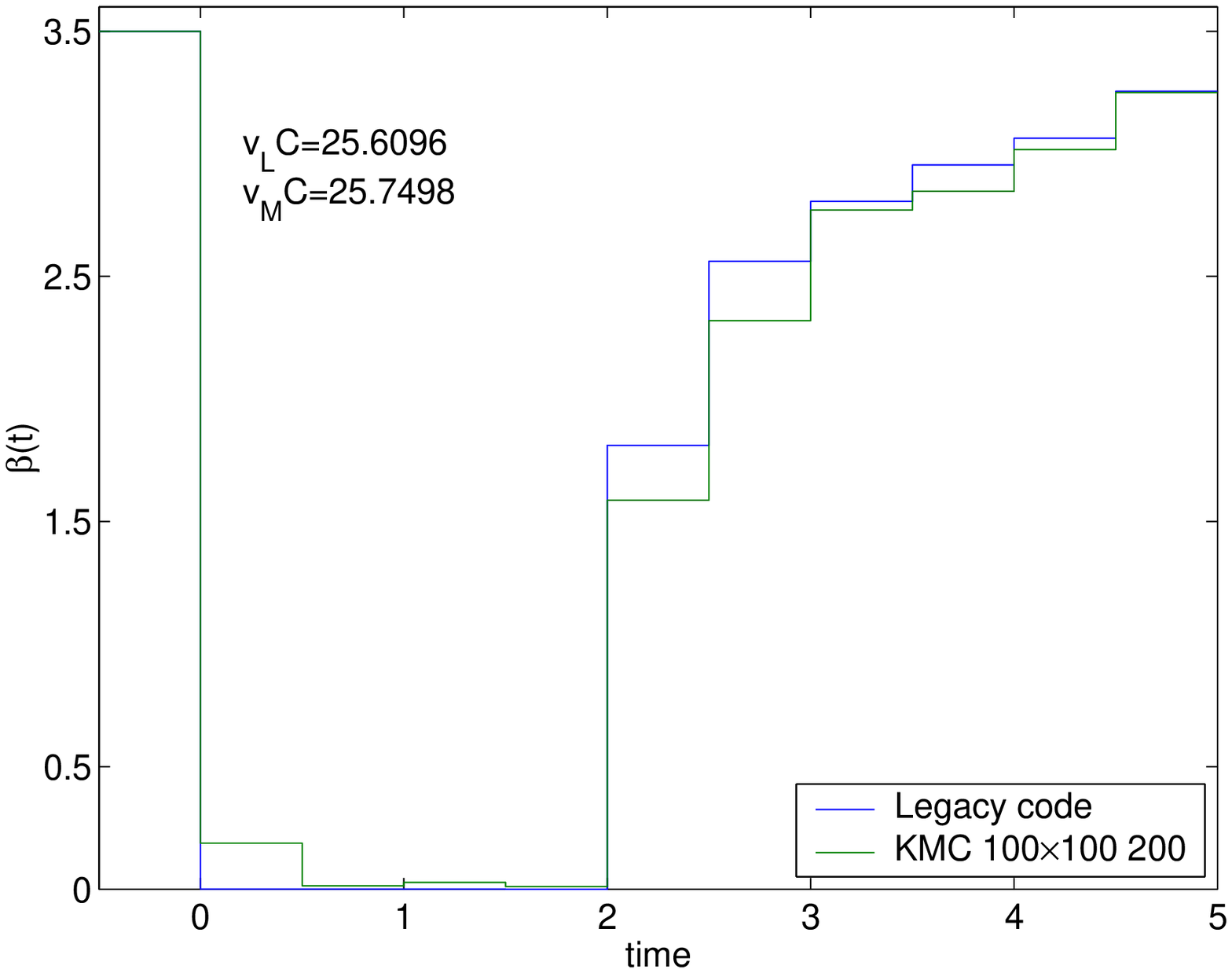,height=2.0in}\hspace{2mm}
            b)\psfig{file=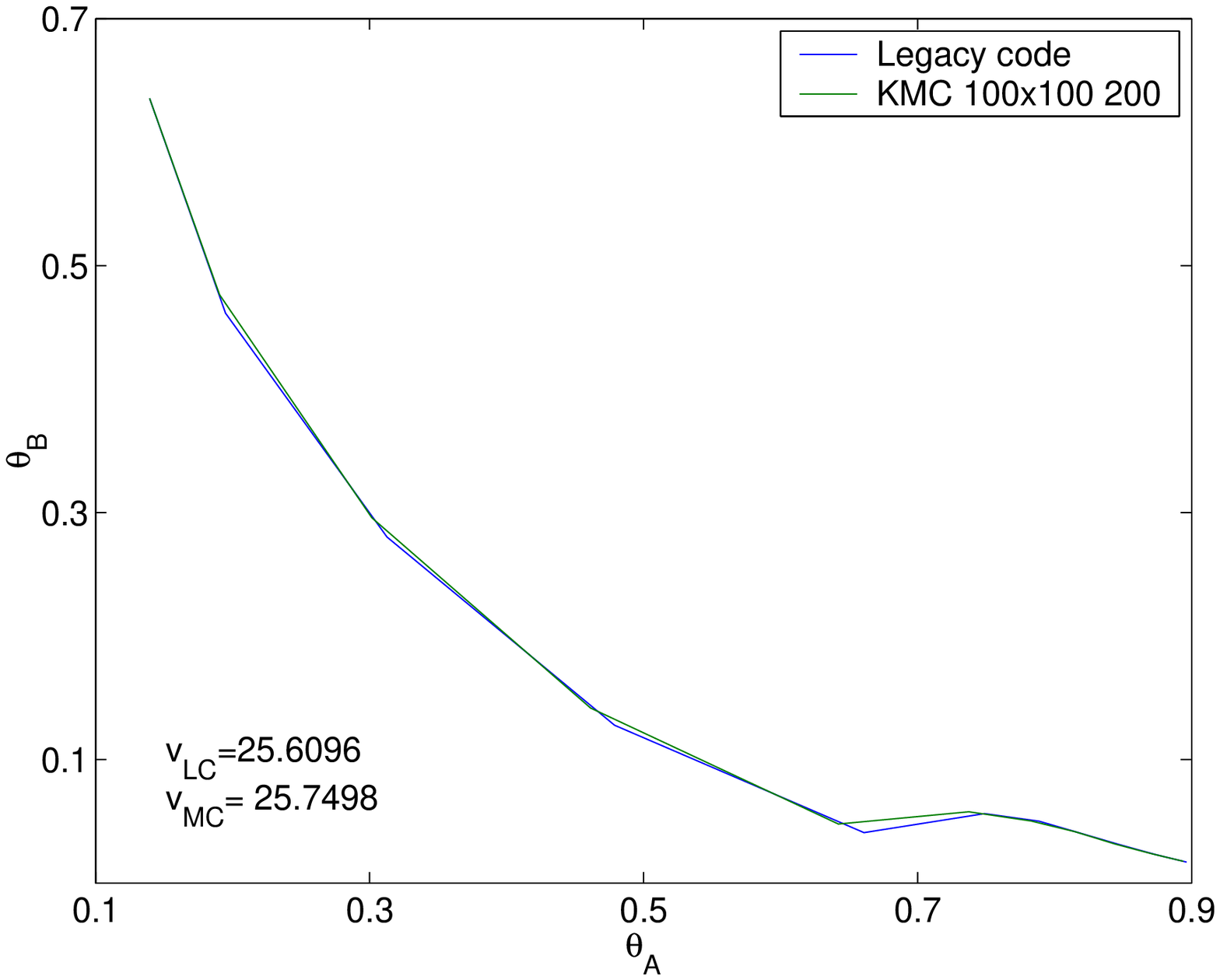,height=2.0in}}%,width=1.2in
\centerline{c)\psfig{file=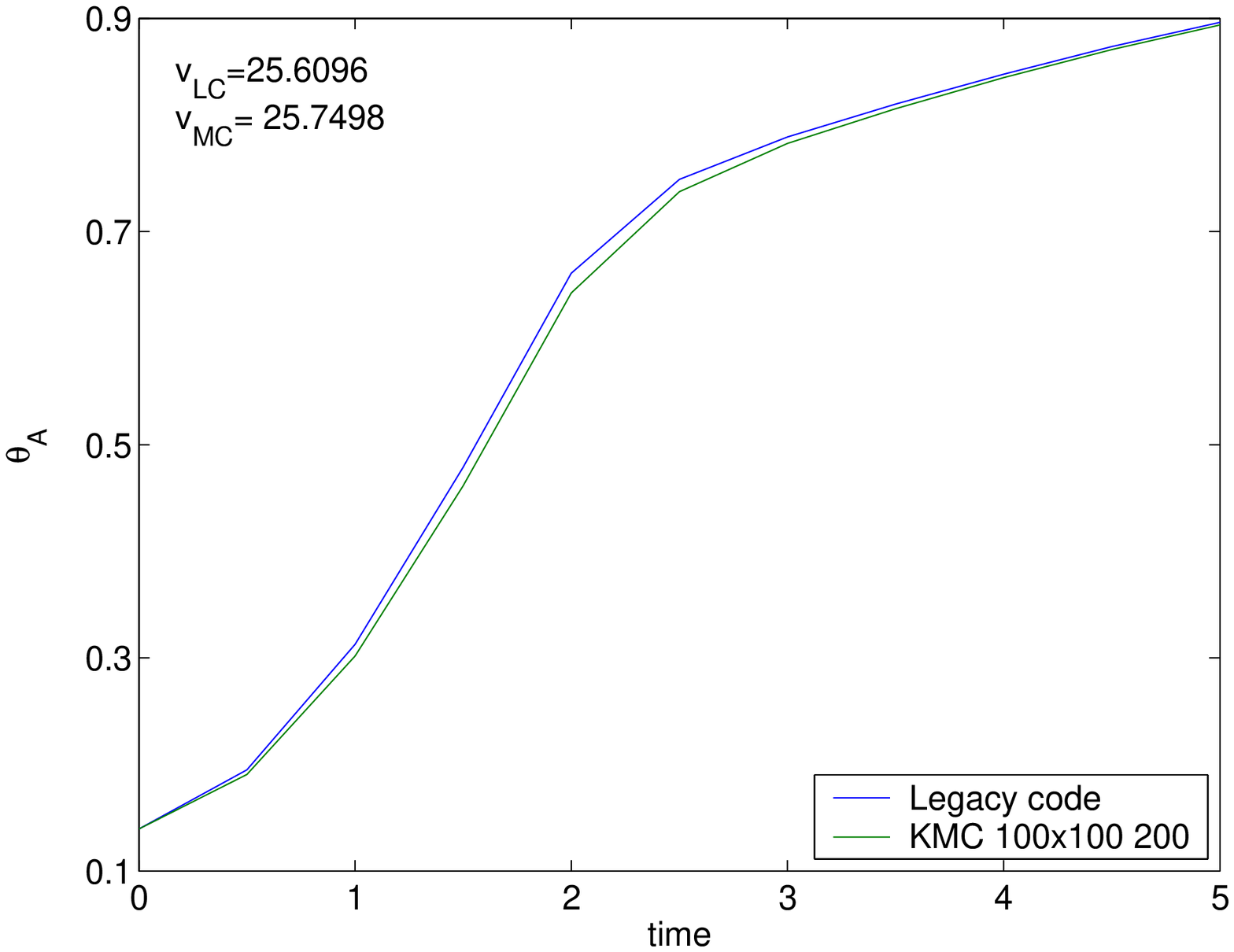,height=2.0in}\hspace{2mm}
            d)\psfig{file=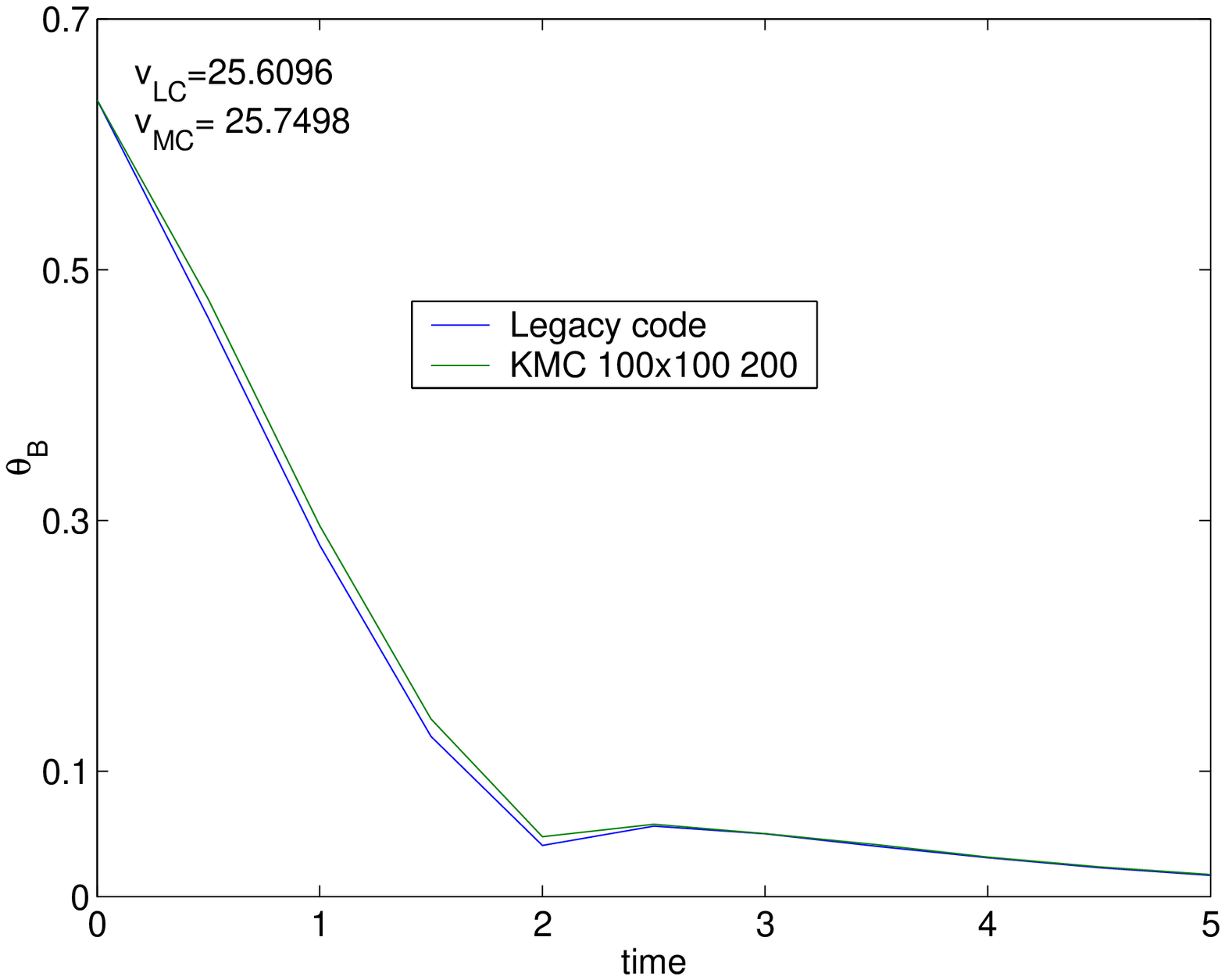,height=2.0in}}%,width=1.2in
\caption{Results of Hooke-Jeeves algorithm through numerical
integration of Eq.\ref{rxn_2ode} (blue line) and using KMC
simulations with $100 \times 100$ lattice size and $200$
repetitions (green line), a) Optimal temporal profile of $O_2$
adsorption rate $\beta$, b) phase portrait of process evolution,
c) evolution of $CO$ coverage $\theta_A$, d) evolution of $O_2$
coverage $\theta_B$ ($t_{f}=5$, $N=10$, $\beta_{ss}=3.5$).}
\label{fig.2D.HJ1}
\end{figure}

The kinetic Monte-Carlo time-stepper and the Hooke-Jeeves
algorithm were used to compute the coarse optimal temporal profile
of the $O_2$ adsorption rate $\beta(t)$ for large timesteps of
$T=0.5$ and lattice size of $100 \times 100$ and averaging of 200
runs to describe the process evolution.
{ For comparison purposes the optimal profile was also identified,
using successive quadratic programming (SQP), when the constraints
were computed through the integration of the mean field
deterministic process model of Eq.\ref{rxn_2ode} (called the
deterministic profile and program, respectively, for the rest of
the section).}
The identified temporal profile of $\beta$ is shown in Figure
\ref{fig.2D.HJ1}a with a blue line and is compared to the computed
profile of the deterministic program, shown with a green line.
The value of the objective function is
$v_{MC,100\times 100, 200}=25.7498$, a relative error of $0.55\%$
over the optimal value of $v_{LC}=25.6096$ . The total time for
the computation of the optimal profile using KMC (first initial
run) was 8 minutes and 16 seconds.
Note that the values of the objective presented
-for comparison purposes, and only upon convergence- are computed by
applying the located optimal policy to the mean field
ODE model of Eq.\ref{rxn_2ode}.
Similar results were also reached when
using the deterministic mean field equations Eq.\ref{rxn_2ode} with
Nelder-Mead after three search initializations and Hooke-Jeeves
also after three search initializations.
Presenting the $CO$ and $O_2$ coverage in Figures
\ref{fig.2D.HJ1}c and \ref{fig.2D.HJ1}d respectively, we observe
that the $CO$ oxidation evolves similarly for both $\beta$
profiles and both realizations (the deterministic and the
KMC-based one).
The same observation is made in
Figure \ref{fig.2D.HJ1}b, the phase portrait of the process
evolution.
We observe that the two evolution profiles, the coarse-grained one
using the KMC computed optimal policy (green line) and the one computed
by integrating Eq.\ref{rxn_2ode} lie close to each other at all times.

To obtain a better resolved approximation of the continuous-time optimal
manipulation of parameter $\beta$, we decreased the time-step to
$T=0.1$.
The discrete-time formulation was subsequently solved using the
Hooke-Jeeves algorithm; results are shown in Table
\ref{t:CO.HJ.MC0.1-obj}.
Specifically, the search algorithm was reinitialized three times,
using the previously estimated optimal profile to initialize the
search, before an optimal over all scales was declared in the case
of the noise-free legacy simulator.
Similarly, when using the KMC
realization of the process, the search was terminated when the
resulting optimal profile was a slight perturbation to the
previous result and the average value of the objective (over ten
simulations) was within error bounds of the previous value.

\begin{figure}[htbp]
\centerline{a)\psfig{file=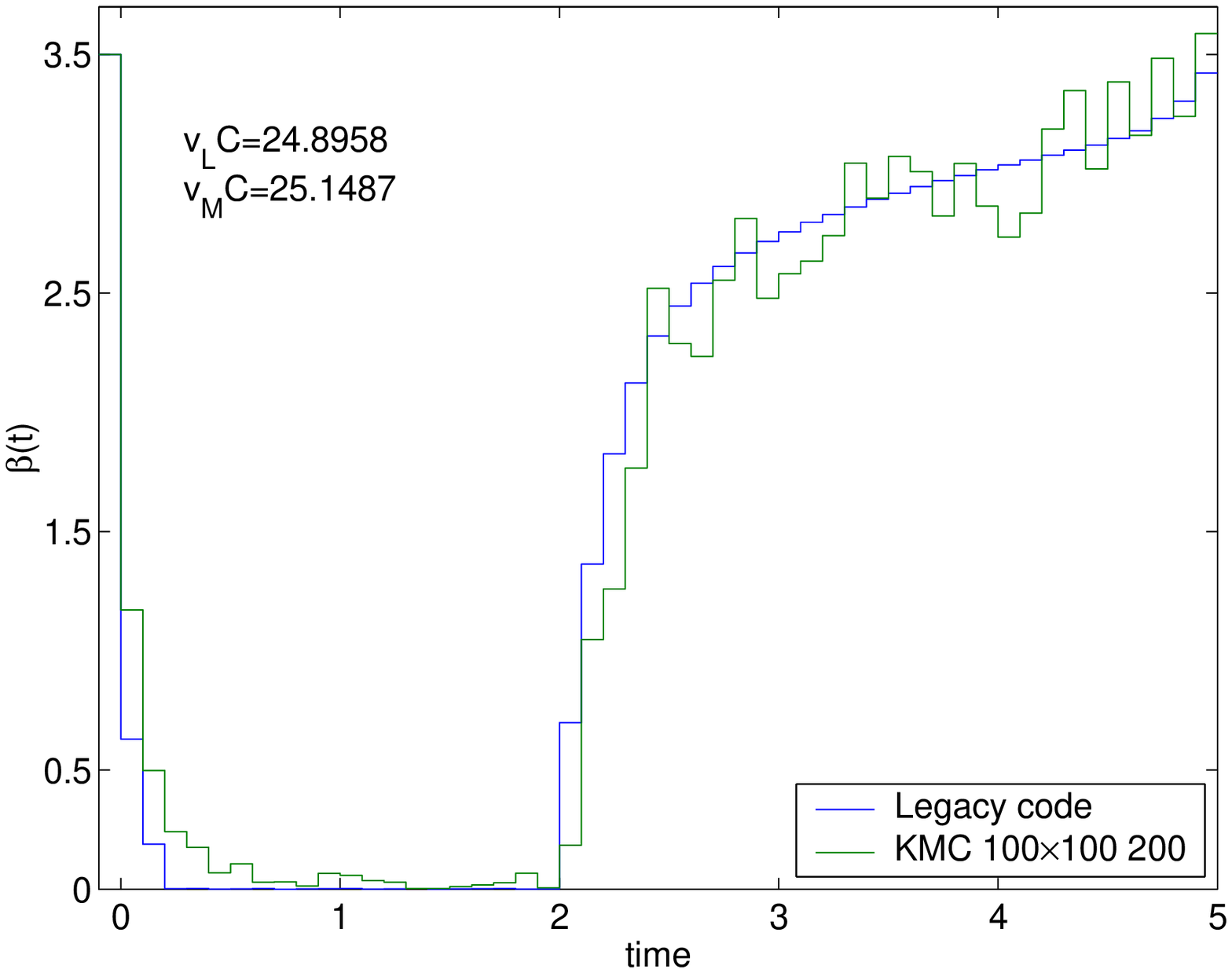,height=2.0in}\hspace{2mm}
            b)\psfig{file=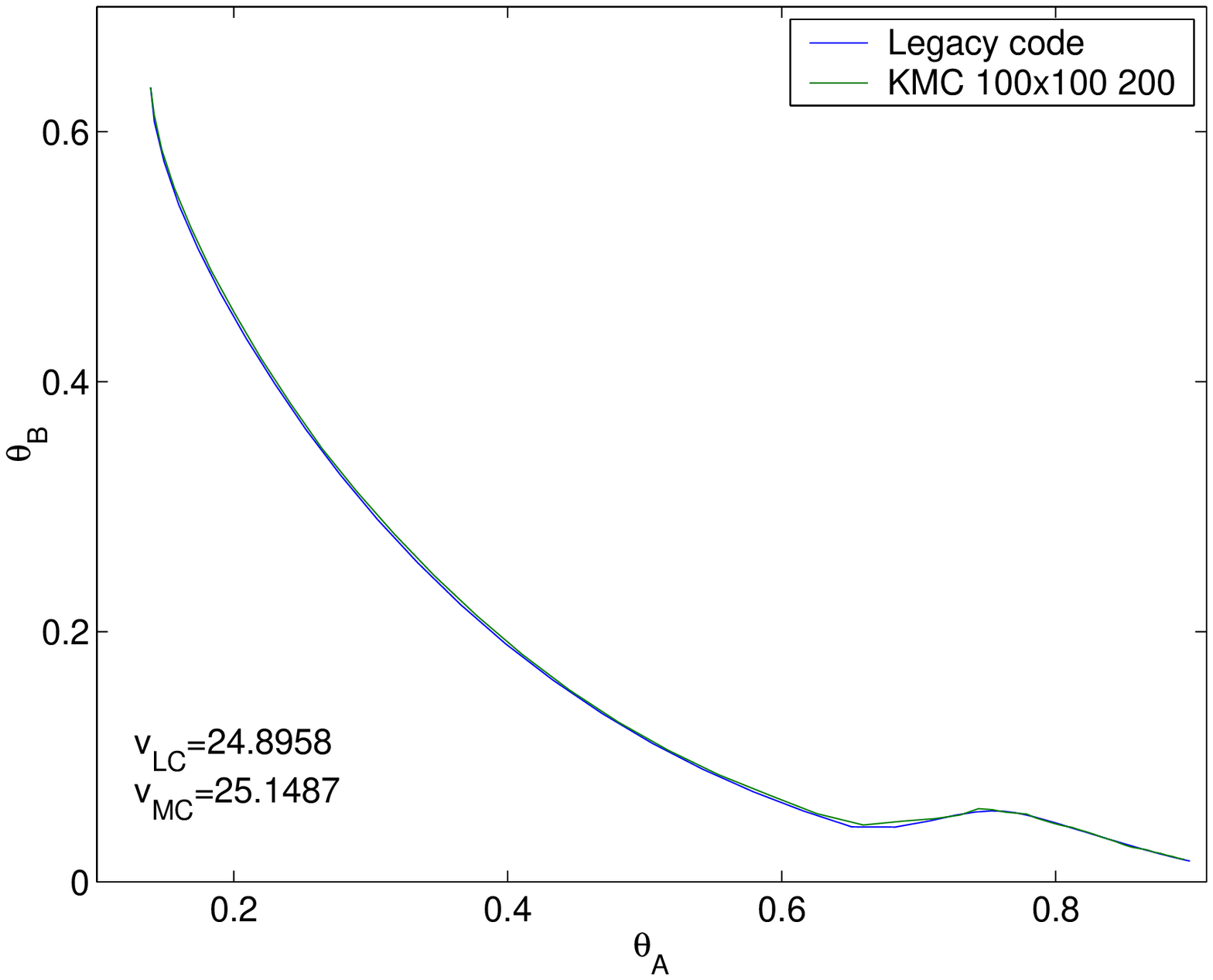,height=2.0in}}%,width=1.2in
\centerline{c)\psfig{file=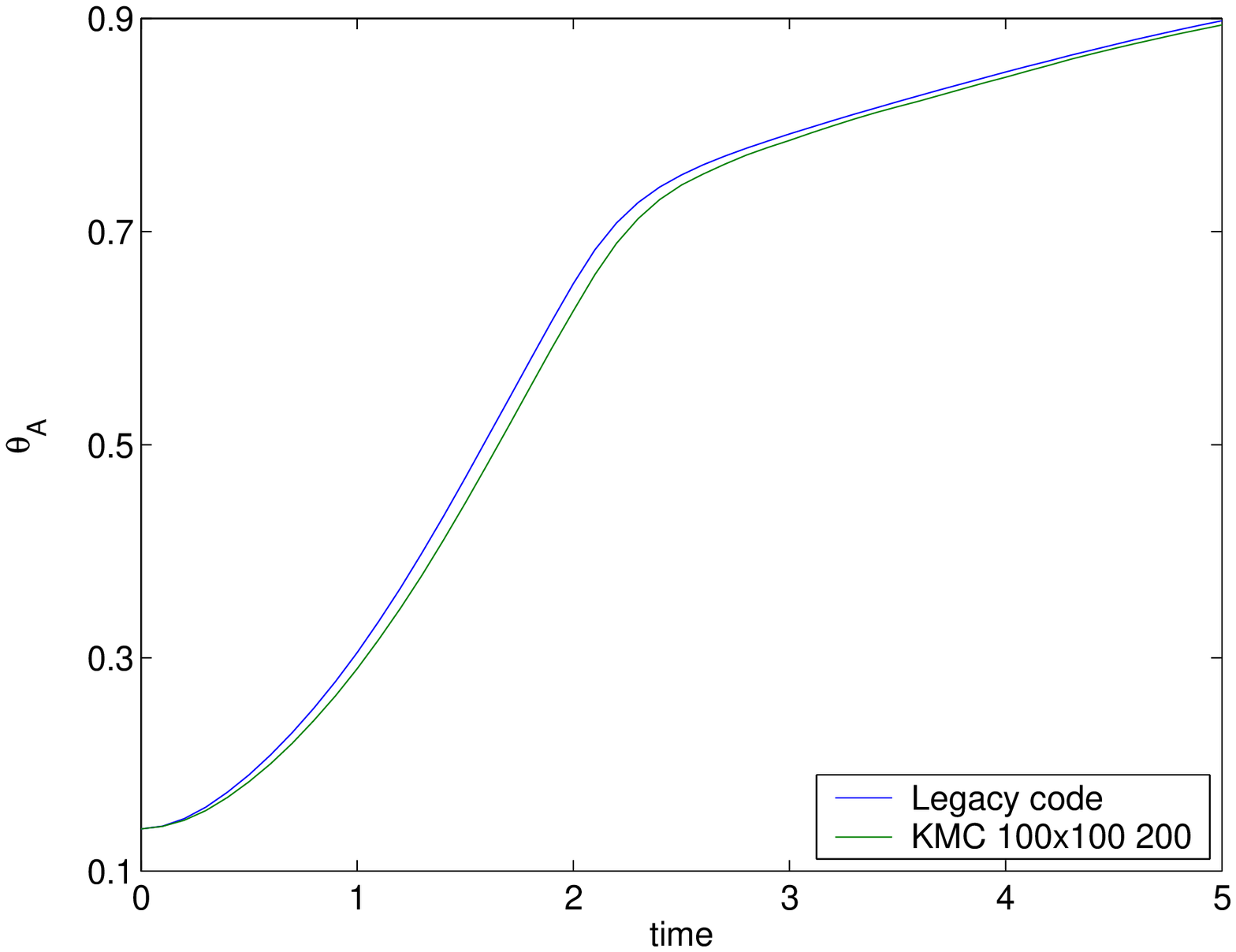,height=2.0in}\hspace{2mm}
            d)\psfig{file=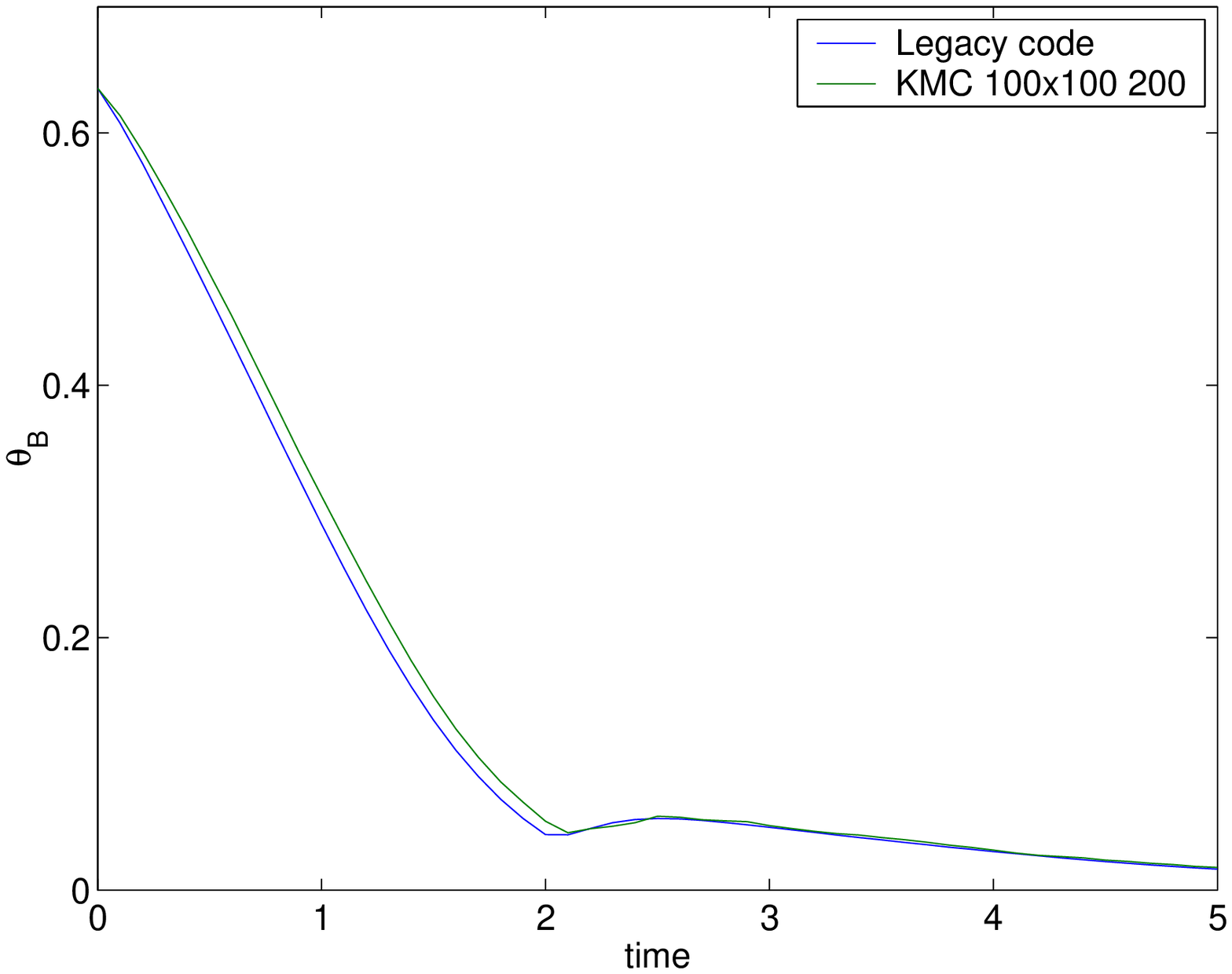,height=2.0in}}%,width=1.2in
\caption{Results of Hooke-Jeeves algorithm through numerical
integration of Eq.\ref{process} (blue line) and using KMC
simulations with $100 \times 100$ lattice size and $200$
repetitions (green line), a) Optimal temporal profile of $O_2$
adsorption rate $\beta$, b) phase portrait of process evolution,
c) evolution of $CO$ coverage $\theta_A$, d) evolution of $O_2$
coverage $\theta_B$ ($t_{f}=5$, $N=50$, $\beta_{ss}=3.5$).}
\label{fig.2D.HJ.MC100@0.1}
\end{figure}

Initially, to solve the optimization problem of
Eq.\ref{opt_fin-cvp} with the Hooke-Jeeves search algorithm we
used a KMC realization of the process with a $N_l=100\times 100$
and $M_r=200$ for the coarse time-stepper output; this combination
identified the $\beta$ temporal profile presented in Figure
\ref{fig.2D.HJ.MC100@0.1}a (green line).
We observe that the
profile is a perturbation around to the resulting $\beta$ profile
(blue line) of the deterministic SQP search.

In Figure \ref{fig.2D.HJ.MC100@0.1}b we present the phase portrait
of the process evolution for the coarse KMC-based HJ computed
profile (blue line); the system response is closely related to the
system behavior under the deterministic, SQP-computed, profile
(green line).
The same observation is also made based on the evolution of $CO$
and $O_2$ coverage, shown in Figures \ref{fig.2D.HJ.MC100@0.1}c
and \ref{fig.2D.HJ.MC100@0.1}d respectively.
%-
The closeness of the two phase portraits under the two different
switching policies suggests that the objective function is
relatively insensitive to the change of $\beta$, a further
difficulty towards the computation of a truly optimal policy.

To better illustrate the effect of noise originating from the
stochastic nature on the time stepper, we run ten independent KMC
simulations with $N_l=100\times 100$, $M_r=200$ for the switching
policy of Figure \ref{fig.2D.HJ.MC100@0.1}.
The average value of the objective function for these ten runs was
$v_{av,MC}=25.1272$ with standard deviation $\sigma_{MC}=0.0241$,
while the value given by the mean-field ODE model is given at
Table \ref{t:CO.HJ.MC0.1-obj}, 1st line.
Based on the intuition gained from these runs, we
subsequently repeated the third pass of the specific search ten
times.
The objective function was then independently computed for the
{\it previously} identified optimal profiles using the ODE model
(for comparison purposes); the average value was computed to be
$v_{av,LC}=25.1928$ with standard deviation $\sigma_{LC}=0.0681$.
Let us reiterate that the KMC coarse time-stepper and the ODE
time-stepper represent {\it different} approximate realizations of
the same process and as such the predictions of the two
time-steppers may lie very close, but in general they are not
necessarily the same.

\begin{figure}[htbp]
\centerline{a)\psfig{file=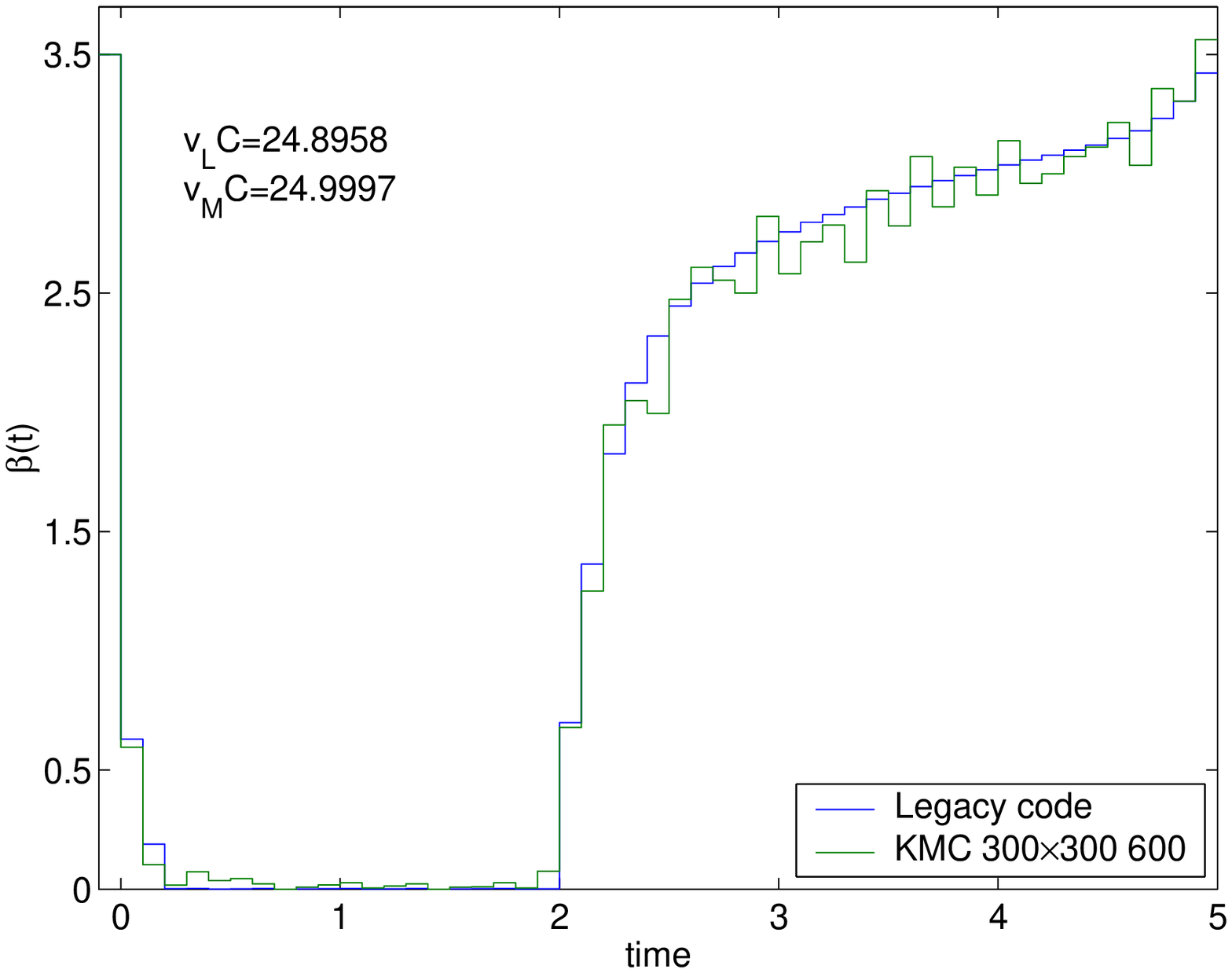,height=2.0in}\hspace{2mm}
            b)\psfig{file=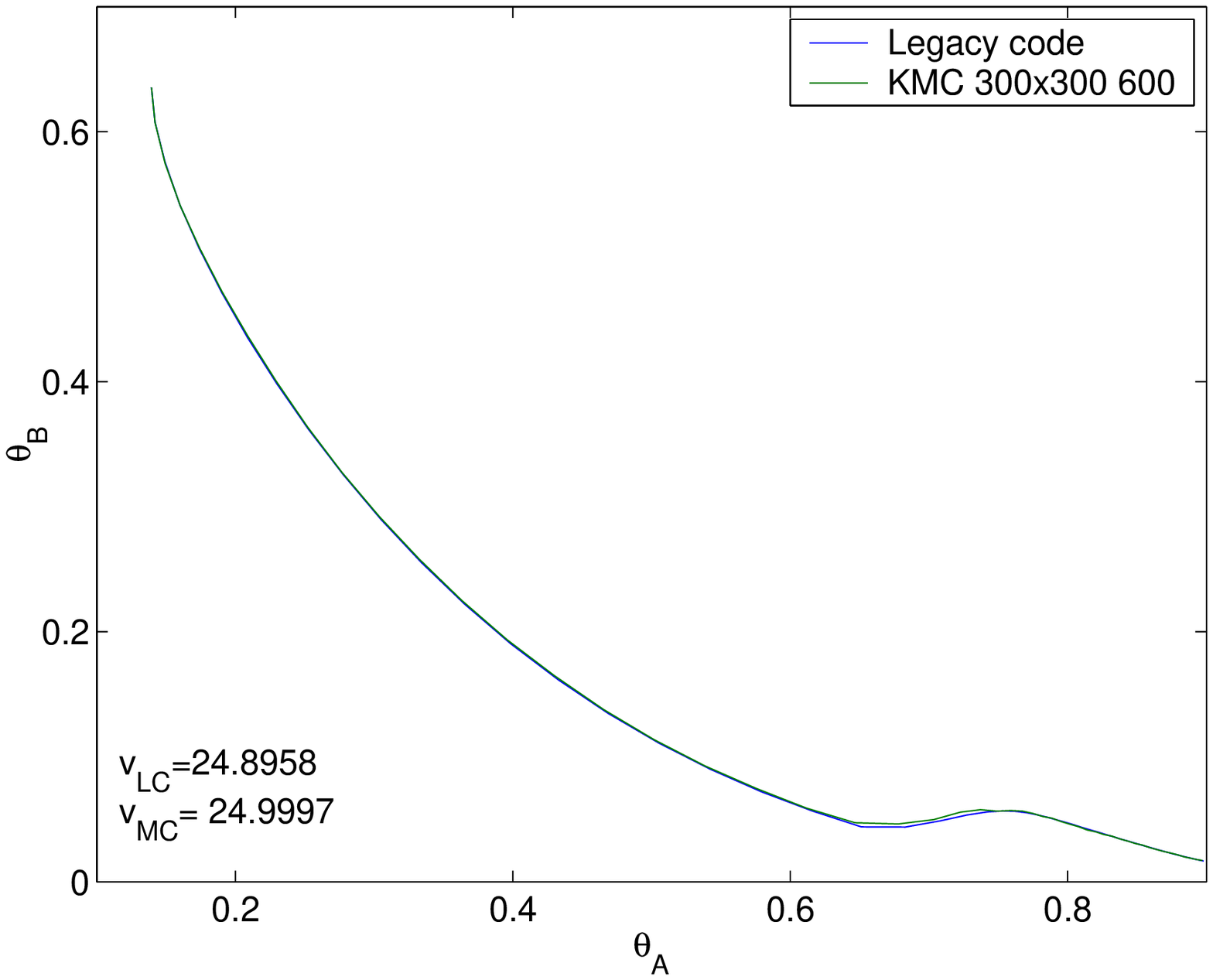,height=2.0in}}%,width=1.2in
\centerline{c)\psfig{file=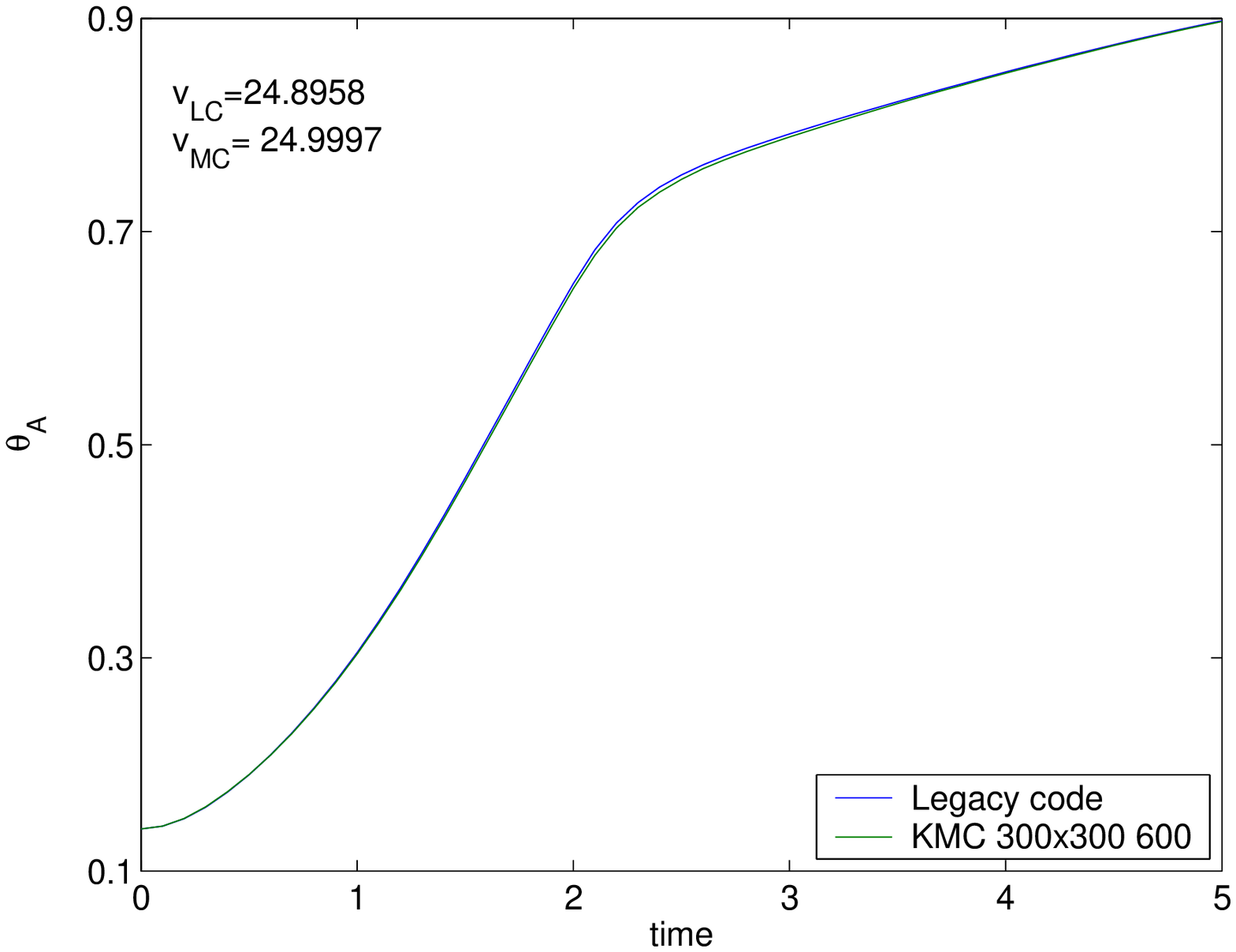,height=2.0in}\hspace{2mm}
            d)\psfig{file=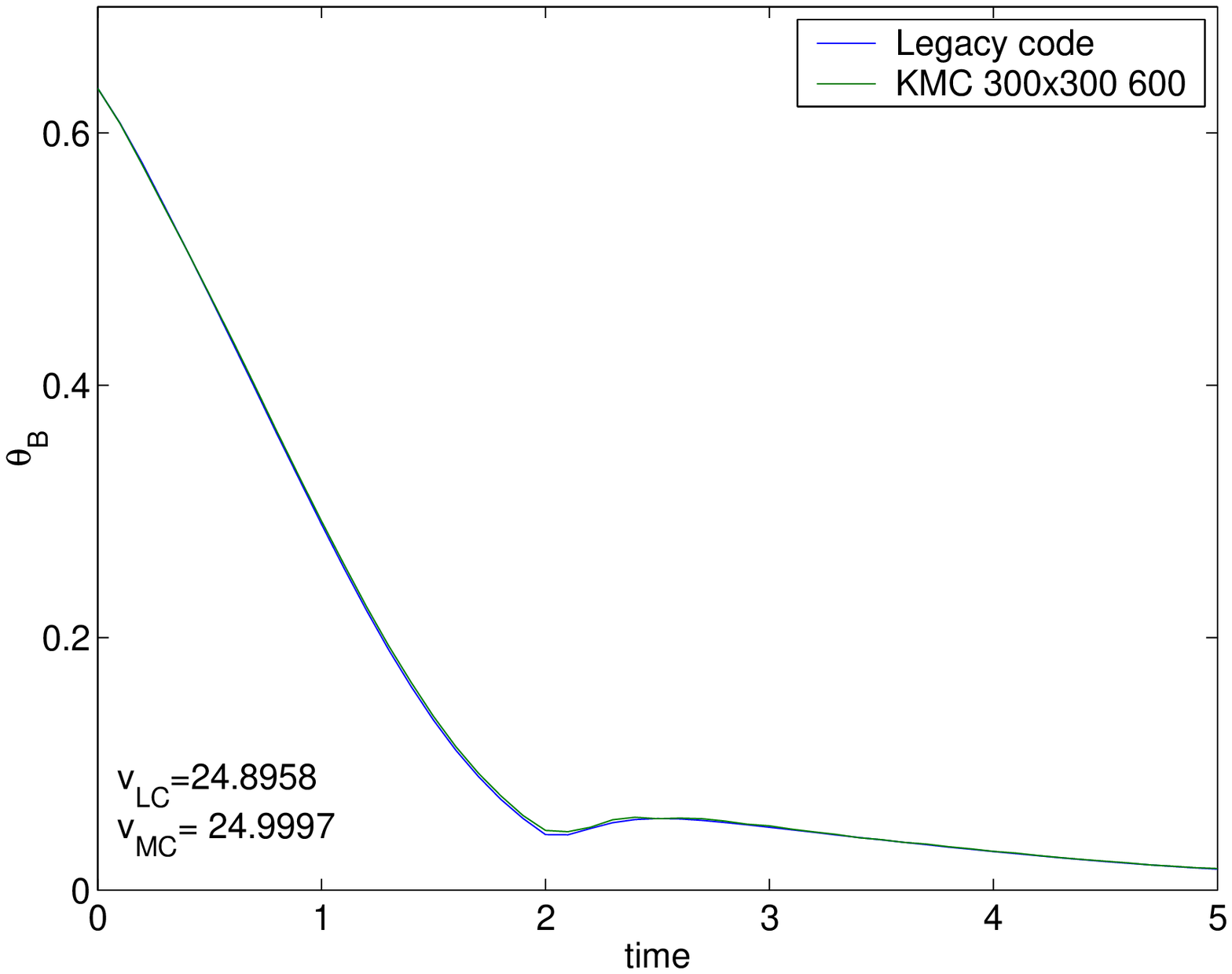,height=2.0in}}%,width=1.2in
\caption{Results of numerical integration of Eq.\ref{rxn_2ode}
(blue line) and  KMC simulations with $300 \times 300$ lattice
size averaged over $600$ runs (green line), a) Optimal temporal
profile of $O_2$ adsorption rate $\beta$, identified by
Hooke-Jeeves search algorithm, b) phase portrait of process
evolution, c) evolution of $CO$ coverage $\theta_A$, d) evolution
of $O_2$ coverage $\theta_B$ ($t_{f}=5$, $N=50$,
$\beta_{ss}=3.5$).} \label{fig.2D.HJ.MC300@0.1}
\end{figure}

The increase in accuracy comes at a high computational cost, as
can be seen from the total CPU time needed for each search, when
using a KMC realization of the process with a $100\times100$
lattice size and an average of $200$ runs, compared to when using
a KMC realization of the process with a $300\times300$ lattice
size and an average of $600$ runs.
In order to reduce the total
computational cost, during the last run (line 3 in Table
\ref{t:CO.HJ.MC0.1-obj}), we initially started a search with the
accuracy of the coarse timestepper implicitly set at a low value
(KMC realization with $100\times100$ lattice size and an average
of $200$ runs).
Once the search converged, we indirectly increased
the accuracy of the timestepper (KMC realization with
$300\times300$ lattice size and an average of $600$ runs) and
initiated a second pass.
The computational savings amounted to
approximately 55000 seconds.
Computational costs can be further reduced by
adjusting the scales during each new search, taking into account the
variance of the coarse time-stepper results ensemble.
Using a KMC realization of the process with a $300\times300$
lattice size and an average of $600$ runs for the coarse
time-stepper output, the Hooke-Jeeves search algorithm identified the
$\beta$ temporal profile presented in Figure
\ref{fig.2D.HJ.MC300@0.1}a (green line).
This compares well with the
$\beta$ profile (blue line) resulting from a deterministic SQP search.

\begin{figure}[htbp]
\centerline{a)\psfig{file=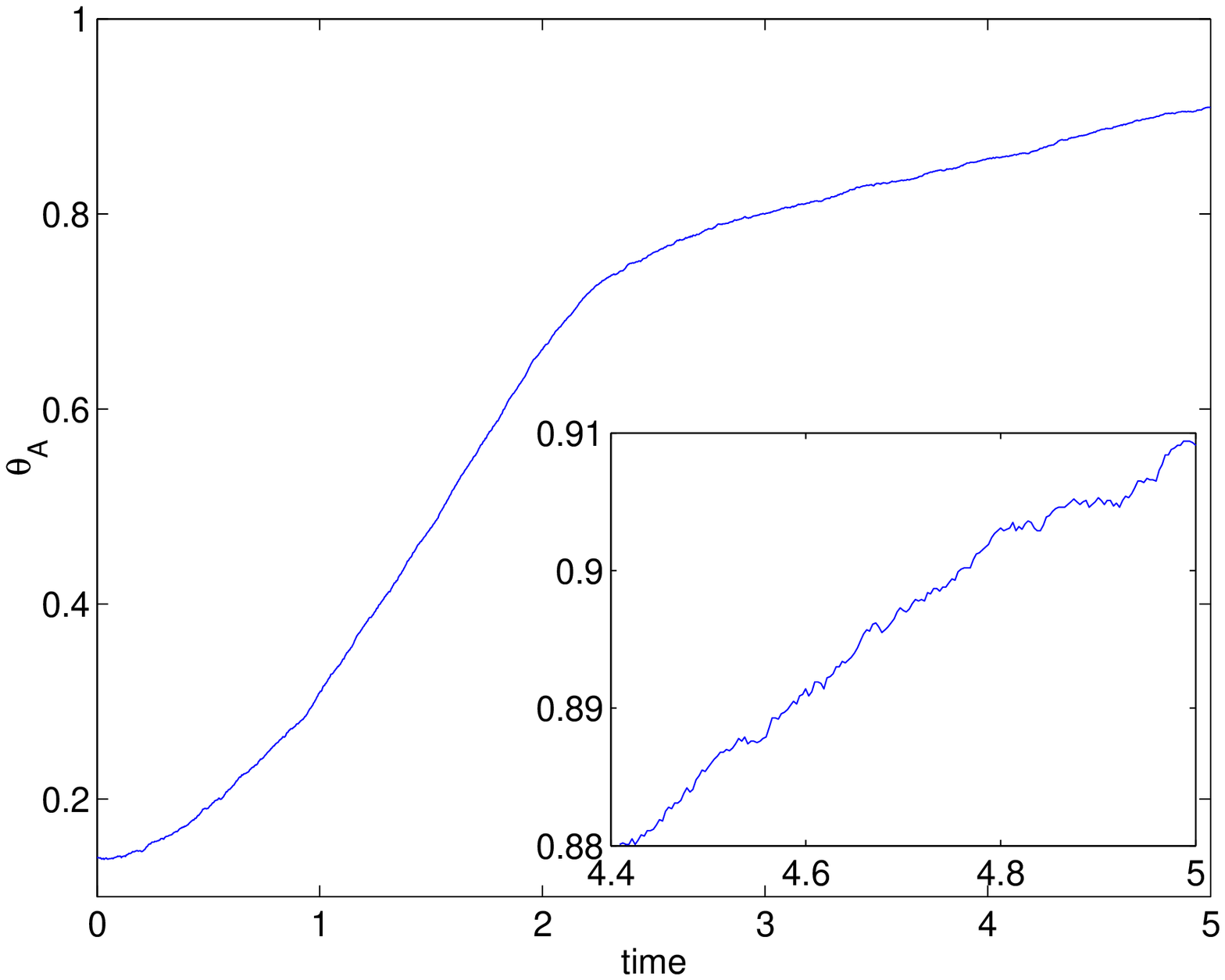,height=2.0in}\hspace{2mm}
            b)\psfig{file=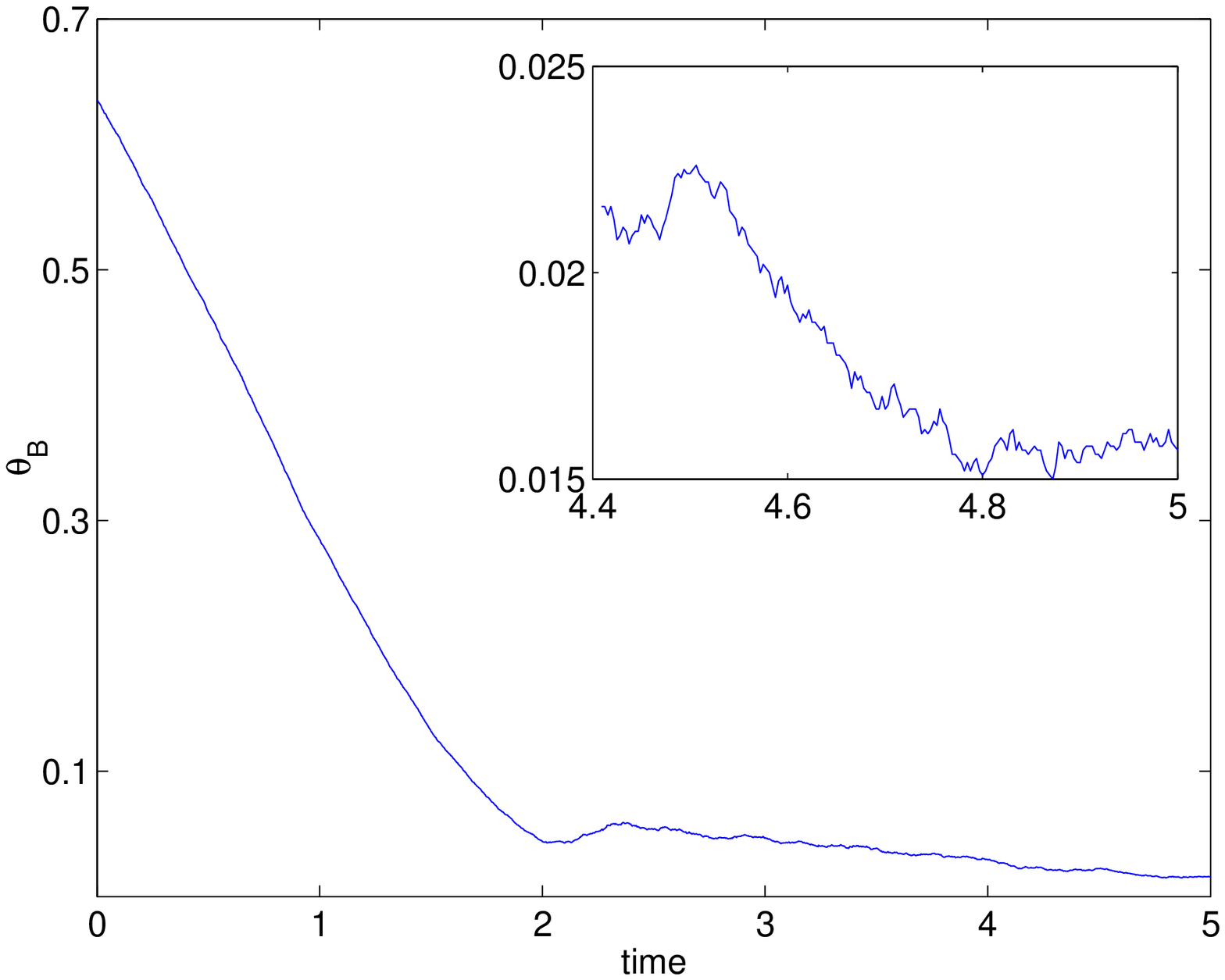,height=2.0in}}%,width=1.2in
\caption{ a) evolution of $CO$ coverage $\theta_A$, b) evolution
of $O_2$ coverage $\theta_B$, for a single KMC realization,
$N_l=100\times 100$, $\delta t=0.0031$, for the Hooke-Jeeves
computed $\beta$ temporal profile ($\beta_{ss}=3.5$).}
\label{fig.2D.HJ.MC300@0.1MC}
\end{figure}

Furthermore, in Figure \ref{fig.2D.HJ.MC300@0.1}b we present the
phase portrait of the process evolution for the coarse-KMC based,
HJ computed profile (blue line) and the deterministic SQP-computed
profile (green line); we observe that the system response is
practically the same for both computed profiles of $\beta(t)$.
%
%%%YANNISNOW
%  I took this out, since it was just said above
%The closeness of the two phase portraits under the two different
%switching policies implies that the objective function is
%relatively insensitive to the change of $\beta$, a further
%difficulty towards the computation of the optimal policy.
%
The same observation is also made based on the evolution of $CO$
and $O_2$ coverage, shown in Figures \ref{fig.2D.HJ.MC300@0.1}c
and \ref{fig.2D.HJ.MC300@0.1}d respectively, for an averaged
realization of the system, and in Figures
\ref{fig.2D.HJ.MC300@0.1MC}a and \ref{fig.2D.HJ.MC300@0.1MC}b
respectively for a single KMC realization, with time-step $\delta
t=0.0031$.

\begin{figure}[htbp]
\centerline{a)\psfig{file=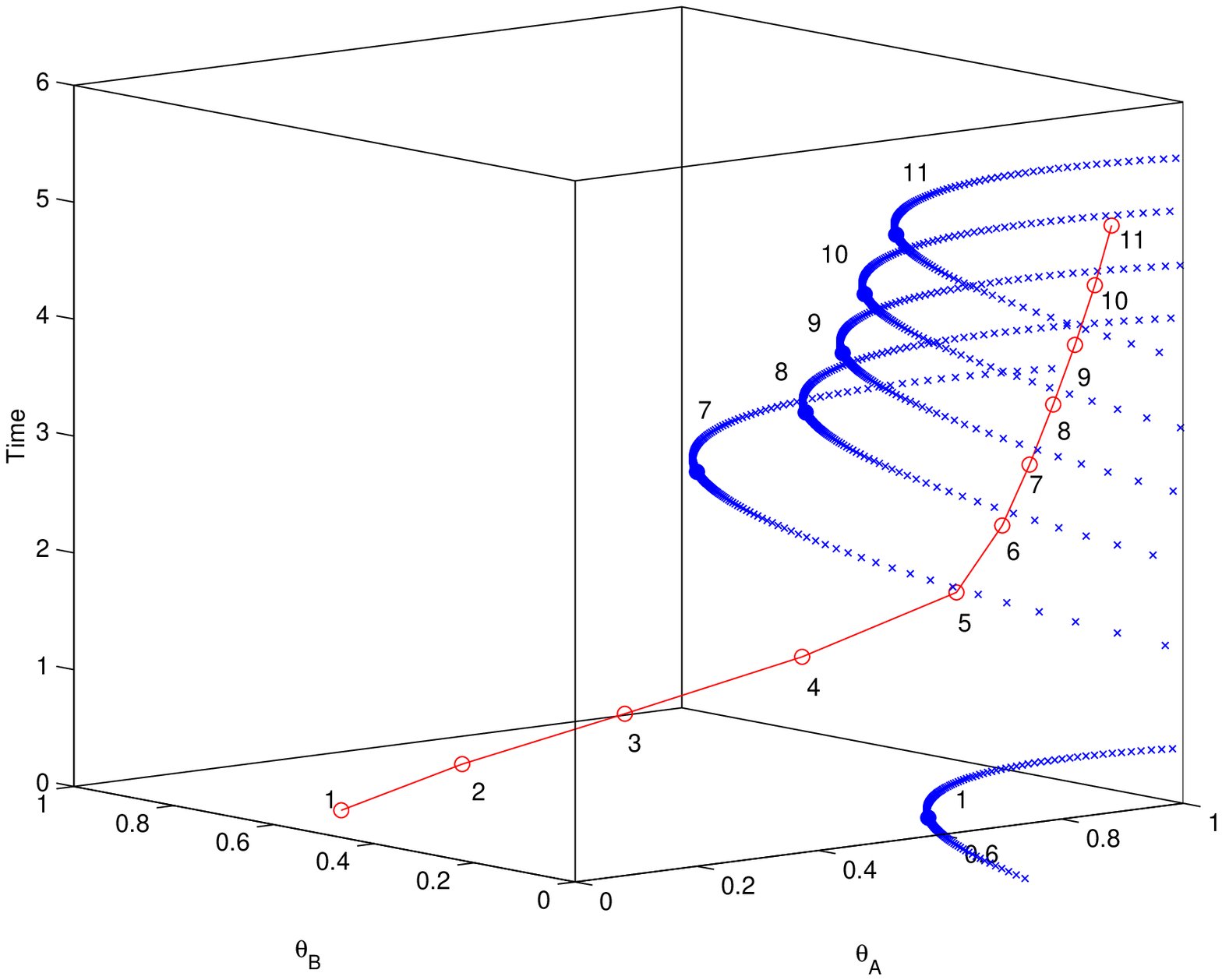,height=2.0in}\hspace{2mm}
            b)\psfig{file=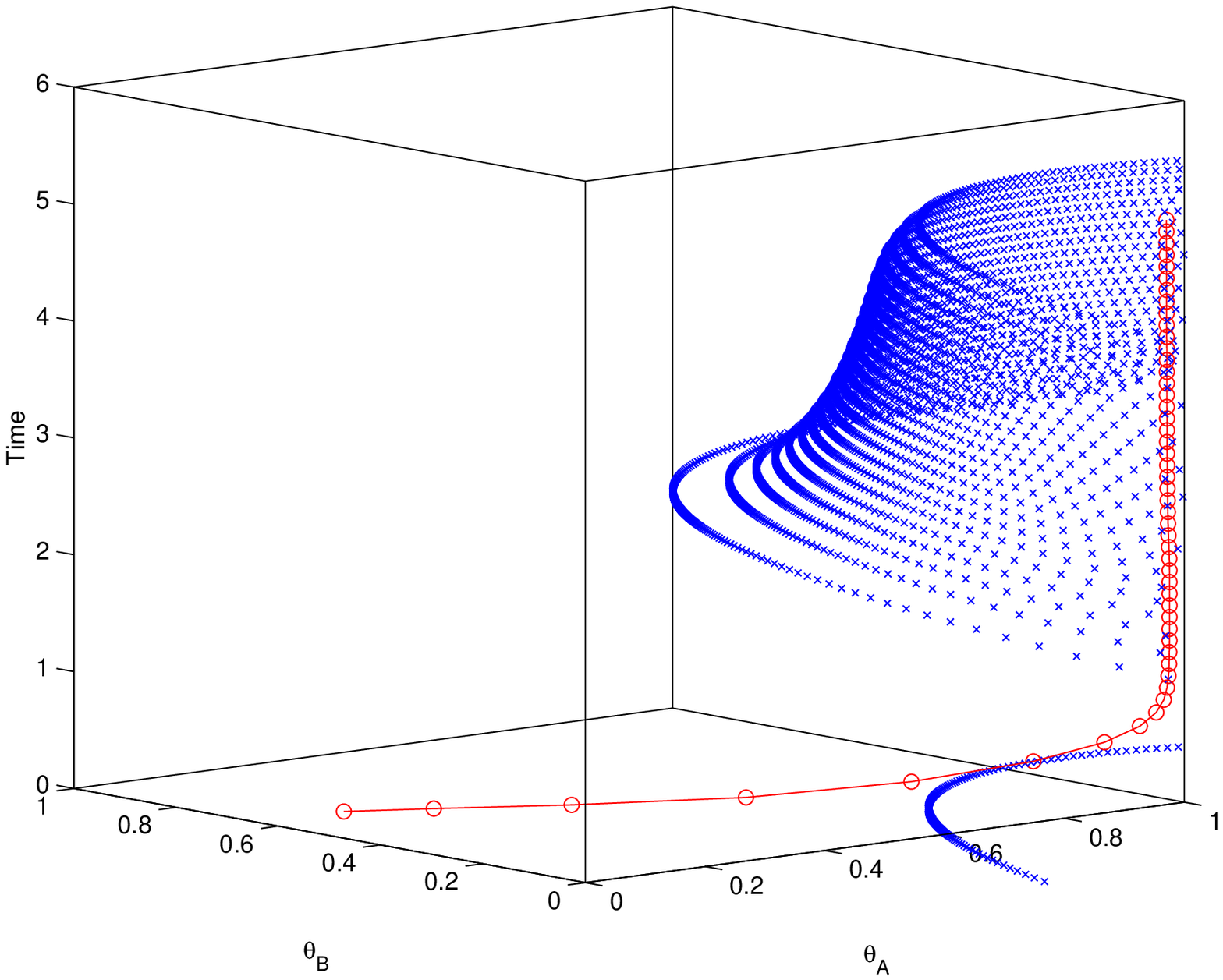,height=2.0in}}%,width=1.2in
\caption{Optimal trajectory for system of Eq.\ref{rxn_2ode} (red
line), identified by Hook-Jeeves search algorithm, and evolution
of separatrix (blue line), a) for $N=10$, b) FOR $N=50$,
identified by ($t_{f}=5$, $\beta_{ss}=3.5$).} \label{fig.3D-CO.HJ}
\end{figure}

The located optimal system trajectory is also shown in Figure
\ref{fig.3D-CO.HJ}a for $T=0.5$ and Figure \ref{fig.3D-CO.HJ}b for
$T=0.1$, with a red line; the blue line represents the separatrix
of the system (if it exists) for the current $\beta(t)$.
The optimal switching policy takes the (expected) phase point from
$\theta_{ss,s}$ through the separatrix of the saddle (unstable) steady state
$\theta_{ss,i}$; note that during traversing the separatrix line
of $\beta_ss$, beta(t) is reduced to a value where no unstable
steady state (and thus separatrix) exists.
After this is accomplished, the optimal $\beta(t)$ rapidly
converges back to $\beta_{ss}$ (see Figure \ref{figIF1}a).
The coarse phase point is now within the region of attraction of
the steady state  $\theta_{ss,f}$, and no particular action is
needed.

Another way to accelerate the computation of the optimal switching
policy, is by solving a sequence of optimization programs,
reducing the time-step in each iteration from an initial
``coarse'' to a final ``finer'' one.
In effect the suggestion is to use a method inspired from
multigrid optimization methods (e.g., \cite{bbdm01a,bbdm01b}).
During each iteration, we utilize the
result of the previous iteration as an initial guess, which we
adapt to the larger number of decisions in time that we need to
take.
Specifically for this work, we used cubic splines to fit the
optimal parameter profiles obtained from $n$-th iteration, and
resampled in the new temporal domain to obtain the initial guess
for the refined search of the $n+1$th iteration.
Subsequently, during the $n+1$th iteration we use an off-the-shelf
search algorithm to obtain optimal process parameter profiles.
As is expected in all multigrid methods, the accuracy of the
off-the-shelf search algorithm used has to be increased to
converge closer to the optimal solution.
This is done in our
case by gradually reducing the size of the scales.
Once we have reached the ``finer'' time-step and obtained the
final solution using this method, one must investigate if the
solution is ``optimal over all scales''.
This is accomplished in
our case by performing a new search using scales with a wide range
of values.
In Table \ref{t:CO.HJ.MCVar-obj} we present the
computational efforts needed, following the aforementioned
procedure, using HJ algorithm with a $100\times 100$ lattice with
$200$ run averaging, KMC realization of the process.
In this case
the scales were kept constant. We observe that we obtain a first
pass of the HJ algorithm at $T=0.1$ in a total of 3983 seconds.

%------------
\section{Conclusions}\label{s:concs}
We presented a computational methodology for the location of
coarse-grained optimal operating parameter policies for chemically
reacting systems described by microscopic/stochastic evolution
rules.
In particular, we approximated optimal operating policies
switching bistable reacting systems from one stationary
state to another; our approach is intended
for systems for which macroscopic, coarse
evolution equations exist but are not available in closed form.
The advantage of the proposed method lies in the establishment,
through the coarse time-stepper, of a computational bridge between
atomistic/ stochastic simulators and traditional (in particular,
derivative free) optimization algorithms.
The approach can be directly extended to systems with higher
dimensional expected behavior (see for example \cite{mmk02}), and
possibly, through matrix-free methods, to systems with infinite
dimensional (spatially distributed, but dissipative) expected
behavior \cite{GKT02}.
Current efforts are focused on applying this methodology to
the study of coarse-grained optimal paths
associated with transitions and rare events in computational
chemistry; in this case the coarse timestepper will be estimated
from processing the results of an ``inner" molecular dynamics or an
``inner" Monte Carlo simulator (\cite{HK02,Dima2}).
In this case, the knowledge of appropriate macroscopic observables
(reaction coordinates) is crucial, and data processing techniques
capable of extracting such reduced data descriptions (e.g.,
\cite{kernelPCA}) become important.
The approach we described in this paper can be easily combined
with ``empirical" observables suggested by such data processing
algorithms.

\section*{Acknowledgements}
Financial support from the Air Force Office of Scientific Research
(Dynamics and Control), National Science Foundation, ITR, and
Pennsylvania State University, Chemical Engineering department, is
gratefully acknowledged. The authors are indebted to professor
Alexei G. Makeev of Moscow State University, department of
Computational Mathematics and Cybernetics for providing the
kinetic Monte Carlo time-steppers used in this work and Ee-Sunn
Chia for the visualization of the CO oxidation trajectories.

%\vspace*{-2mm}
\bibliographystyle{plain}
\bibliography{refOP.MC,refOP.Opt}

%%%%----------------Tables----------------------------------------------

\clearpage
\begin{table}[h]
\caption{NO reduction process parameters}\label{t:proc} %\vspace{1mm}
\begin{center}
\begin{tabular}{|c|c|cc|}
\hline
 Parameter  & Value & \multicolumn{2}{|c|}{Steady states}\\ \hline
 $T$        & 0.25 & $\theta_{ss,s}$ & 0.3301 \\
 $N$        & 20   & $\theta_{ss,i}$ & 0.6803 \\
 $t_f$      & 5    & $\theta_{ss,f}$ & 0.9896 \\
 $k_{ss}$   & 0.45 &&\\
 $\alpha$   & 1.0  &&\\
 $\gamma$   & 0.01 &&\\
\hline
 \end{tabular}
\end{center}
\end{table}

\clearpage
\begin{table}[h]
\caption{Hooke-Jeeves search results}\label{t:HJ.MC-obj}%\vspace{1mm}
\begin{center}
\begin{tabular}{l}
\begin{tabular}{|c||c|c|c|c|}
\hline
 Model  & $N_l$ & $M_r$ & Objective & t$^*$ [s]\\
\hline\hline
 Legacy & $-$              & $-$   & 10.3709 & 129   \\
 KMC    & $100 \times 100$ & $200$ & 10.5418 & 674   \\
 KMC    & $200 \times 200$ & $400$ & 10.4155 & 4673  \\
 KMC    & $300 \times 300$ & $600$ & 10.3973 & 16300 \\
 KMC    & $400 \times 400$ & $800$ & 10.3815 & 38469 \\
 KMC    & $500 \times 500$ &$1000$ & 10.3811 & 75307 \\
 \hline
\end{tabular}\\
$^*$ single CPU pentium $IV$ at $2.4\; GHz$
\end{tabular}\end{center}
\end{table}

\clearpage
\begin{table}[h]
\caption{Implicit Filtering search results}\label{t:IF.MC-obj}%\vspace*{0.5mm}
\begin{center}
\begin{tabular}{l}
\begin{tabular}{|c||c|c|c|c|}
\hline
 Model  & $d_{max}$   & Scales                & Objective & t$^*$ [s]\\%$^{\dagger}$
\hline\hline
 Legacy & $-$         & $-$                   & 10.3709   & 168  \\
 KMC    & $0.005$     & $2^{-0},\dots,2^{-3}$ & 10.5461   & 297  \\
 KMC    & $0.005$     & $2^{-0},\dots,2^{-4}$ & 10.7922   & 282  \\
 KMC    & $0.0005$    & $2^{-3},\dots,2^{-5}$ & 10.4110   & 3684 \\
 KMC    & $0.001$     & $2^{-0},\dots,2^{-4}$ & 10.4129   & 6858 \\
 \hline
\end{tabular}\\
 $^*$ single CPU pentium $IV$ at $2.0\; GHz$
\end{tabular}\end{center}
\end{table}

\clearpage
\begin{table}[h]
\caption{CO oxidation process parameters}\label{t:COproc} %\vspace{1mm}
\begin{center}
\begin{tabular}{|c|c|cc|}
\hline
 Parameter   & Value & \multicolumn{2}{|c|}{Steady states}\\ \hline
 $T$         & 0.25 & $\theta_{CO,ss,s}$  & .13944 \\
 $N$         & 20   & $\theta_{CO,ss,i}$  & .67526 \\
 $t_f$       & 5    & $\theta_{CO,ss,f}$  & .97101 \\
 $k_{r}$     & 1.0  & $\theta_{O_2,ss,s}$ & .63553 \\
 $\alpha$    & 1.6  & $\theta_{O_2,ss,i}$ & .11452 \\
 $\gamma$    & 0.04 & $\theta_{O_2,ss,f}$ & .00137 \\
 $\beta_{ss}$& 3.5 &&\\
\hline
 \end{tabular}
\end{center}
\end{table}

\clearpage
\begin{table}[h]
 \caption{Hooke Jeeves search results. $T=0.1$, $t_f=5$}\label{t:CO.HJ.MC0.1-obj}
\begin{center}
\begin{tabular}{l}
\begin{tabular}{|c||c|c|c|c|c|}
\hline
 Model  & $N_l$       & $M_r$    & Passes   & Objective & Total t [s]\\%$^{\dagger}$
\hline\hline
 Legacy & $-$           & $-$        & 3        & 24.8958   & 3381  \\
 KMC    & $100\times100$& $200$ runs & 3        & 25.1487   & 8182  \\
 KMC    & $300\times300$& $600$ runs & 2        & 24.9997   & 57482$^*$  \\
 \hline
\end{tabular}\\
 $^*$First search used a $100\times100$ lattice size KMC realization.
 \end{tabular}\end{center}
\end{table}

\clearpage
\begin{table}[h]
 \caption{Hooke Jeeves search results. $t_f=5$, KMC with $100\times100$
  lattice, 200 runs average}\label{t:CO.HJ.MCVar-obj}
\begin{center}
\begin{tabular}{l}
\begin{tabular}{|c|c|c|c|c|}
\hline
 Model  & Timestep  & Objective & Total t [s]\\%$^{\dagger}$
\hline\hline
 KMC    & $0.50$    & 25.7498   &  481 \\
 KMC    & $0.25$    & 25.4800   &  899 \\
 KMC    & $0.10$    & 25.2734   & 2603 \\
 KMC    & $0.10^*$  & 25.2607   & 3100 \\
 \hline
\end{tabular}\vspace{1mm}\\
% $^*$First search used a $100\times100$ lattice size KMC realization.
$^*$ Second pass at desired timestep.
\end{tabular}\end{center}
\end{table}

%%%%----------------Figures---------------------------------------------

\clearpage

\clearpage

\clearpage

\clearpage

\clearpage

\clearpage

\clearpage

\clearpage

\clearpage

\clearpage

\clearpage

\end{document}